\documentclass[10pt]{article}

\usepackage{amsthm}
\usepackage{enumitem,multirow,supertabular,tabularx ,algorithm,inputenc,comment,booktabs,makecell,multicol}
\usepackage{caption,parskip}
\usepackage[dvipsnames]{xcolor}
\usepackage[numbers]{natbib}
\usepackage{array}
\usepackage{xspace}
\usepackage{setspace}
\definecolor{mycolor}{rgb}{.5,.5,.5}
\usepackage[hidelinks,colorlinks = true,linkcolor = BrickRed,urlcolor  = blue, citecolor = blue, anchorcolor = blue]{hyperref}

\usepackage[margin=1in]{geometry}

\usepackage{graphicx,algorithm,multirow,booktabs}
\usepackage[]{subcaption}
\usepackage[dvipsnames]{xcolor}
\usepackage{amsmath,amsfonts,amssymb,enumitem}

\DeclareMathOperator*{\argmax}{arg\,max}

\newcommand{\beql}[1]{\begin{equation}\label{#1}}

\newcommand{\eeql}{\end{equation}}
\newcommand{\eqn}[1]{(\ref{#1})}

\newtheorem{theorem}{Theorem}
\newtheorem{assumption}[theorem]{Assumption}

\newtheorem{lemma}[theorem]{Lemma}
\newtheorem{proposition}[theorem]{Proposition}
\newtheorem{remark}[theorem]{Remark}
\theoremstyle{remark}

\newcolumntype{L}[1]{>{\raggedright\let\newline\\\arraybackslash\hspace{0pt}}m{#1}}
\newcolumntype{C}[1]{>{\centering\let\newline\\\arraybackslash\hspace{0pt}}m{#1}}
\newcolumntype{R}[1]{>{\raggedleft\let\newline\\\arraybackslash\hspace{0pt}}m{#1}}
\usepackage{titlesec}
\titlespacing*{\section}
{0pt}{1.5ex plus 1ex minus .2ex}{.7ex plus .2ex}
\titlespacing*{\subsection}
{0pt}{1.5ex plus 1ex minus .2ex}{.7ex plus .2ex}
\title{Reward Maximization in General Dynamic Matching Systems} % CAN WE CHANGE THE TITLE

\author{Mohammadreza Nazari
        \thanks{ISE Department, Lehigh University, Bethlehem, PA, USA.  E-mail: {\tt mon314@lehigh.edu}}
        \and
        Alexander L. Stolyar
        \thanks{ISE and CSL Department, University of Illinois at Urbana-Champaign, Urbana, IL, USA.  E-mail: {\tt stolyar@illinois.edu}}}
\newcommand*{\titleGP}{\begingroup % Create the command for including the title page in the document
\centering % Center all text
\vspace*{\baselineskip} % White space at the top of the page

{\LARGE%\scshape % Small caps
Reward Maximization in General Dynamic Matching Systems\par} % Location and year
\vspace*{1.5\baselineskip} % Whitespace between location/year and editors
\begin{table}[H]

\footnotesize
  \centering
    \begin{tabular}{C{7cm} C{7cm} }
{\Large Mohammadreza Nazari\par} &{\Large Alexander L. Stolyar\par}\\
 Lehigh University&University of Illinois at Urbana-Champaign\\
200 West Packer Ave.& 1308 W. Main Street, 156CSL\\
Bethlehem, PA 18015&Urbana, IL 61801\\
mon314@lehigh.edu &stolyar@illinois.edu 
    \end{tabular}%
\end{table}%

\endgroup}

\begin{document}

\titleGP 
\begin{abstract}
We consider a matching system with random arrivals of items of different types. The items wait in queues -- one per each item type -- until they are ``matched.'' Each matching requires certain quantities of items of different types; after a matching is activated, the associated items leave the system. There exists a finite set of possible matchings, each producing a certain amount of ``reward''. This model has a broad range of important applications, including assemble-to-order systems, Internet advertising, matching web portals, etc.

We propose an optimal matching scheme in the sense that it asymptotically maximizes the long-term average matching reward, while keeping the queues stable. The scheme makes matching decisions in a specially constructed virtual system, which in turn control decisions in the physical system. The key feature of the virtual system is that, unlike the physical one, it allows the queues to become negative. The matchings in the virtual system are controlled by an extended version of the greedy primal-dual (GPD) algorithm, which we prove to be asymptotically optimal -- this in turn implies the asymptotic optimality of the entire scheme.
The scheme is real-time, at any time it uses simple rules based on the current state of virtual and physical queues. It is very robust in that it does not require any knowledge of the item arrival rates, and automatically adapts to changing rates.

The extended GPD algorithm and its asymptotic optimality apply to a quite general queueing network framework, not limited to matching problems, and therefore is of independent interest.
\end{abstract}

\section{Introduction}\label{sec.introduction}
We consider a dynamic matching system with random arrivals. Items of different types arrive in the system according to a stochastic process and wait in their dedicated queues to be ``matched.'' 
Each matching requires certain quantities of items of different types; after a matching is activated, the associated items leave the system. There exists a finite number of possible matchings, each producing a certain amount of ``reward''.  The objective is to maximize long-term average rewards, subject to the constraint that the queues of currently unmatched items remain stochastically stable. In this paper we propose a dynamic matching scheme and prove its asymptotic optimality. (In fact, the policy works for a more general objective, being a concave function of the long-term rates at which different matchings are used.)

\begin{figure}[h]
\centering
	\includegraphics[width=0.4\textwidth]{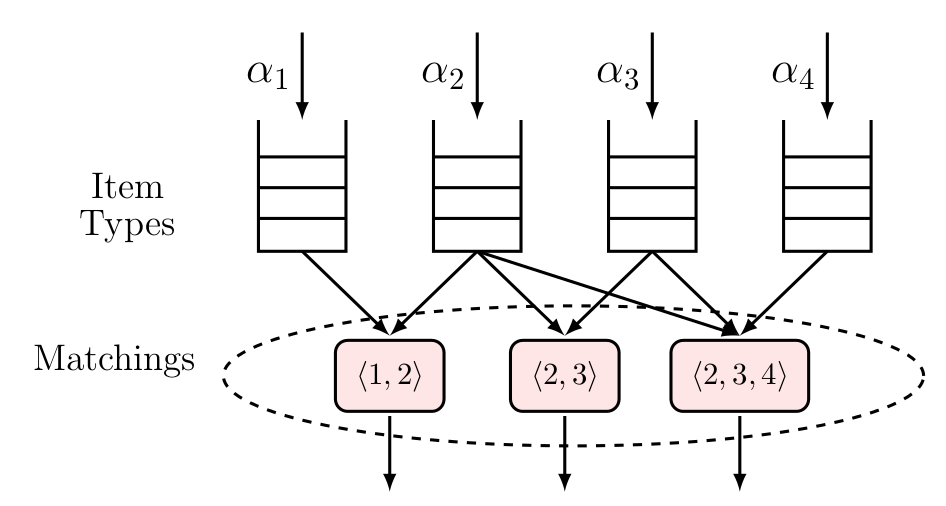}	
\caption{An example of the matching model}
\label{figmatch}
\end{figure}

Figure~\ref{figmatch} shows an example of a matching system with 4 item types. The items 
arrive as a random process, as individual items or in batches. The average arrival rate of type $i$ items is
$\alpha_i$. There exist 3 possible matchings; e.g. $\langle 1,2\rangle$ is a matching which matches one item of type 1 with one item of type 2. $\langle2,3,4\rangle$ is another matching which matches one item of types 2, 3 and 4. (In general, unlike in this example, a matching may require more than one item of any given type.) A matching can only be applied if all contributing items are present in the system; and if it is applied, the contributing items instantaneously leave the system.

%\vspace{-.3cm}

The analysis of \emph{static} matching has a large literature (see, e.g., \cite{lovasz2009matching}). The \emph{dynamic} model, which we focus on, has attracted a lot of attention recently, due to large variety of new (or relatively new) important applications. One example is assemble-to-order systems (see e.g. \cite{PW2008} and references therein), where randomly arriving product orders are ``matched'' with sets of parts required for the product assembly. Another application is to Internet advertising \cite{mehta2012online}, where the problem is to find appropriate matchings between the ad slots and the advertisers. Web portals as places for business and personal interactions is an important application; the problem in these portals (such as dating websites, employment portals, online games) is to match people with similar interests \cite{buke2015stabilizing}. Matching problems also arise in systems with random arrival of customers and servers;
for example, in taxi allocation, where matched ``items'' are passengers and taxis \cite{kashyap1966double}. Further applications also can be found in \cite{busic2014optimization,rene2009fcfs}.

Different control objectives may be of interest for matching systems. Gurvich and Ward \cite{gurvich2014dynamic} study the problem of minimizing finite-horizon cumulative holding costs for a model very close to ours. Plambeck and Ward \cite{PW2008}, in the context of assemble-to-order systems, consider a model where item arrival rates can be controlled via a pricing mechanism; the objective includes queueing holding costs in addition to rewards/costs associated with order fulfillments, parts salvaging and/or expediting. Paper \cite{PW2008}, in particular, proposes and studies a discrete-review policy; it involves solving an optimization problem at each review point.

A special case of the matching system, which received considerable attention, is where customers and servers are randomly arriving in the system and each server can be matched with one customer from a certain subset. This model, also known as the (stochastic) bipartite matching system, was initially studied by \cite{rene2009fcfs}. Majority of the previous research for this model was focused on finding the stationary distribution \cite{adan2015reversibility,adan2012exact} and stability issues \cite{buke2015stabilizing,busic2010stability,mairesse2014stability}. Bu{\v s}i\'c et al. \cite{busic2010stability} established the necessary and sufficient conditions for stabilizability of such systems, and have shown that the well known MaxWeight algorithm achieves maximum stability region. The problem of minimizing the long-term average holding cost for the bipartite matching system is studied by \cite{busic2014optimization}. They have shown that with known arrival rates (and some other conditions on the problem structure), a threshold-type policy is asymptotically optimal in the (appropriately defined) heavy traffic regime.

In this paper, we show that the reward-maximizing optimal control of the matching model can be obtained by putting it into a typical queueing network framework. Our scheme uses a specially constructed virtual system, whose state, along with the state of the physical system, determines control decisions via a simple rule. In the virtual system any matching can be applied at any time and the queues are allowed to be negative. The matchings in the virtual system are controlled by (an extended version of) the Greedy Primal-Dual (GPD) algorithm \cite{stolyar2005maximizing}, which maximizes a queueing network utility subject to stability of the queues. Negative queues in the virtual system can be interpreted as the shortages of physical items of the corresponding types. The GPD algorithm in \cite{stolyar2005maximizing} does not allow negative queues, so it is insufficient for the control of our virtual system. {\em The main theoretical contribution} of this paper is that we introduce and study an {\em extended version of GPD,} 
labeled EGPD, which {\em does allow  negative queues,} and prove its asymptotic optimality 
under non-restrictive conditions that we specify.
The approach of using a virtual system to control the original one has been used before, e.g. in \cite{stolyar2010control}, but the virtual system employed in this paper is substantially different, primarily because it allows negative queues.

Our proposed scheme is very robust in that it does not require a priori knowledge of item arrival rates, and automatically adjusts if/when the arrival rates change. It also covers a wide range of applications and control objectives. For example, in the context of assemble-to-order systems, the objective can include rewards/costs
associated with order fulfillments, parts salvaging and/or expediting. 

Although our scheme is designed (and proved asymptotically optimal) for the reward maximization objective, which does {\em not} include holding costs, we will discuss heuristic approaches to how the scheme can be used to 
achieve good performance in terms of 
a more general objective (including holding costs).

%COM SAY THAT WE COVER REWARD MAX WITHIN SINGLE FRAMEWORK, WHICH IS SIMPLE AND VERY ROBUST. DISCRETE REVIEW HAS SHORTCOMINGS

The paper is organized as follows. Section~\ref{sec-notation} contains notation used throughout the paper.
In Section \ref{secmatchmodel} we formally introduce the matching model and the reward maximization problem;
here we also formally define the corresponding virtual system and the overall control scheme, in which the matching algorithm for the virtual system is a key part. 
In Section \ref{sec3network}, we introduce the Extended Greedy Primal-Dual (EGPD) algorithm for a general network 
model, with queues that may be negative, 
and prove asymptotic optimality of EGPD; here we also show that the virtual system algorithm (in Section~\ref{secmatchmodel}) is a special case of EGPD and thus is asymptotically optimal.
(A reader interested mostly in applications of our proposed scheme may skip Section~\ref{sec3network},
at least at first reading.) We evaluate the performance of our scheme via simulations in Section \ref{sec:experiment}.
Finally, in Section~\ref{profit}, we discuss heuristics on how a more general objective, including holding costs,
 can be addressed by tuning EGPD parameters. Some conclusions are given in Section~\ref{sec-conc}.

%%%%%%%%%%%%%%%%%%%%

\section{Basic Notation}
\label{sec-notation}

We denote by $\mathbb R^{}$, $\mathbb R^{}_+$ and $\mathbb R^{}_-$ the set of real, real non-negative and real non-positive numbers, respectively. $\mathbb R^{N}$, $\mathbb R^{N}_+$ and $\mathbb R^{N}_-$ are the corresponding $N$-dimensional vector spaces. A vector $x\in\mathbb R^{N}$ is often written as $x=(x_n,n\in\mathcal{N})$, where $\mathcal{N} = \{1,2,\cdots, N\}$. For two vectors $x,y\in\mathbb R^{N}$,
$ %\begin{equation}
x\cdot y = \sum_{n=1}^{N} x_n y_n %\nonumber
$ %\end{equation}
 is the scalar (dot) product; vector inequality $x \le y$ is understood component-wise.
 The standard Euclidean norm of $x$ is denoted by $\|x\| = \sqrt{x\cdot x}$. The distance between point $x$ and set $V \subseteq \mathbb R^{N}$ is denoted by $\rho(x,V) = \inf_{y\in V} \|x-y\|$. 

 For a vector function $f:\mathbb R_+\rightarrow \mathbb R^{N}$ and a set $V\subseteq \mathbb R^{N}$, the convergence $f(t)\to V$ means that $\rho(f(t),V)\to 0$ as $t\to \infty$.
 
For differentiable functions $f:\mathbb R\rightarrow\mathbb R$ and $g:\mathbb R^{N}\rightarrow\mathbb R$, we use $f'(t)$ (or $(d\slash dt)f(t)$) to denote the derivative with respect to $t$ and $\nabla g(x) =((\partial\slash \partial x_n)g(x),n\in\mathcal{N})$ is the gradient of $g$ at $x\in\mathbb R^{N}$.

For a set $V$ and a real-valued function $g(v)$, $v\in V$,
 \begin{equation}
 \argmax_{v\in V} g(v)\nonumber
 \end{equation}
denotes the subset of vectors $v\in V$ which maximizes $g(v)$.

For $\xi,\eta\in\mathbb R$ and $\gamma\in\mathbb R_+$, we denote: $\xi \wedge \eta=\min\left\{\xi,\eta\right\}$, $\xi\vee\eta=\max\left\{\xi,\eta\right\}$; $\xi^+ = \xi \vee 0$, $\xi^- = (-\xi) \vee 0$;
$\left[\xi\right]^+_\gamma=\xi$ if $\gamma>0$ and $\left[\xi\right]^+_\gamma=\max\{\xi,0\}$ if $\gamma=0$.

Abbreviation \textit{a.e.} means \emph{almost everywhere} with respect to Lebesgue measure.

%%%%%%%%%%%%%%%%%%%%%%%%%%%%%%%%%%%%%%%%%%%%%%%%%
\section{Optimal Control of the Matching System}\label{secmatchmodel}

The outline of this section is as follows. First, we formally define the \emph{physical} matching system in Section \ref{secmatchmodeldef} and discuss the flexibility of this model to include a large variety of practical systems in 
Section~\ref{sec-model-flexibility}. In Section~\ref{secmatchvirtualdef} we introduce a virtual system, corresponding to the physical one. In Section~\ref{secdynvirtual} we define a control scheme, 
such that a certain algorithm runs on the virtual system, and
control decisions for the physical system depend on those in the virtual one. 
We propose a specific algorithm for the virtual system in Section~\ref{secmatchmodelalg};
this algorithm is asymptotically optimal in the sense that,
under certain non-restrictive conditions,
 when the algorithm parameter ($\beta$) goes to zero,
our entire physical/virtual control scheme maximizes average matching reward in the physical system.
(The asymptotic optimality will be proved later, in Section~\ref{sec3network}.)
We discuss features of the virtual system algorithm, and the conditions for its asymptotic 
optimality in Section~\ref{sec-virtual-alg-discussion}. 

\subsection{Definition of the Physical Matching System}
\label{secmatchmodeldef}

Consider a matching system with $I$ \emph{item types} forming set $\mathcal{I} =\{1,\cdots,I\}$. 
The customers arrive in {\em batches}, consisting of items of same or different types.
To simplify exposition, assume that batches arrive as Poisson process, with each batch type chosen upon arrival, independently, according to some fixed distribution. There is a finite number of possible batch types.
The average rate at which type $i$ customers arrive into the system is $\alpha_i>0$.

There is a finite set $\mathcal{J} =\{1,\cdots,J\}$ of possible \emph{matchings}. 
Let $\mu(j)=\left(\mu_i(j),i\in\mathcal{I}\right)$, where $\mu_i(j) \ge 0$ is the required number of type $i$ items to form matching $j\in\mathcal{J}$. Without loss of generality, we can and do
assume that the ``empty'' matching, with all $\mu_i=0$, is an element of $\mathcal{J}$; the empty matching 
is denoted $\langle \emptyset \rangle$. If a matching requires either zero or one item of each type,
it is denoted by the subset of the required item types; say, $\langle 1,2 \rangle$ denotes the matching requiring
one item of type $1$ and one item of type $2$.

Without loss of generality, we can and do assume that the matching decisions are made only at the times of batch arrivals into the system. Essentially without loss of generality, we also assume that at those times at most $m \ge 1$ matchings can be done. To simplify exposition, we further assume that $m=1$ -- it will be clear from our analysis that all results and (with very minor adjustments) proofs hold for arbitrary fixed $m$. Therefore, from now on 
we consider the system as operating in discrete (slotted) time  $t=0, 1, 2, \ldots$,
with i.i.d. batches arriving at those times, and exactly one (possibly empty) matching activated at each $t$.

Further, without loss of generality, we adopt the convention that the items arrived at time $t$ are only available for matching
at time $t+1$. (If items arriving at time $t$ are immediately available for matching, the convention still holds if we simply pretend that they arrived at time $t-1$, after the matching decision at time $t-1$ was made.)
\iffalse
if arrival times are $\tau_t, t=1,2, \ldots$, 
exactly one matching (possibly empty) is done at the times $\tau_t +$ right after the arrivals.
%These times are indexed by the corresponding $t$. 
From now on we consider the system as operating in discrete time $t=1, 2, \ldots$, 
with the 'state at time $t$' meaning the state at time $\tau_t +$ (right after $\tau_t$) in the actual system;
by convention, the 'arrivals at $t$' refers to arrivals at time $\tau_{t+1}$ in the actual system.
\fi

Type $i\in \mathcal{I}$ items waiting to be matched form a first-come-first-served (FCFS) queue; its length is denoted $\hat Q_i$.
At any time $t$, any one matching $j\in\mathcal{J}$ can be activated subject to the constraint that all the required items must be available in the system. 
With activation of matching $j\in\mathcal{J}$,
\begin{enumerate}[label=(\roman*),font=\itshape]
	\item  Certain (real-valued) reward $w_j$ is generated;
	\item Number $\mu_i(j)$ of items is removed from the queues of the corresponding types $i$.
\end{enumerate}
Let $X_j$ be the long-term average reward generated by matching $j$, under a given control policy. We are interested in finding a dynamic matching policy, which maximizes a continuously differentiable
concave utility function $G(X_1,\cdots,X_J)$ subject to the constraint that all queue lengths $\hat Q_i(\cdot)$ remain stochastically stable. Informally speaking, 
stochastic stability means that as time goes to infinity the queues do not ``run away'' to infinity, i.e. remain $O(1)$. Formally, by stochastic stability we will understand positive recurrence of the underlying Markov process, describing the system evolution. (For example, if the process is a countable-state-space irreducible Markov chain, positive recurrence is equivalent to the existence of unique stationary probability distribution and to ergodicity.)
Therefore, stochastic stability ensures that all arriving items are matched, without the backlogs and waiting times of unmatched items building up to infinity over time.

\begin{remark}
\label{remark:stability} \normalfont
\textbf{Stability and long-term averages.} 
We will give a specific definition of long-term average rewards $X_j$ later. 
When the process is Markov, positive recurrent, then $X_j$ can be thought of 
as the {\em steady-state average reward $u_j$} due to type $j$ matchings -- we will elaborate on the relation 
between $X_j$ and $u_j$ later.
\end{remark}

\begin{remark}
\label{remark:negative-mu} \normalfont
\textbf{More general $\mu_i(j)$.} 
Our model and the results hold -- as is -- in the case when the values of $\mu_i(j)$ can be real numbers of any sign. A negative $\mu_i(j)$ means that matching $j$ {\em adds} $|\mu_i(j)|$ items to $\hat Q_i$, 
and by convention any negative number of items of any type is always available for matching completion.
We assume in this paper that $\mu_i(j)$ are non-negative integers to keep the exposition intuitive.
\end{remark}

\subsection{Model Flexibility}
\label{sec-model-flexibility}

The matching model defined in Section~\ref{secmatchmodeldef} is very flexible to include a variety of systems and their features. Let us consider assemble-to-order systems as an example. In such systems, orders for multiple products arrive as a random process. Each product requires a certain number of components of each type to be assembled. Components also arrive into the system as a random process. A product can only be assembled when all necessary parts are available; in which case it brings a certain reward (profit). This is a matching system where the components and product-orders of different types are ``items'', a completed product is a matching comprising one corresponding product-order and the required number of parts. Salvaging and/or disposing of the components is easily accommodated; namely, salvaging/disposing of one component,
labeled as a type $i$ item, can be treated as a matching $\langle i \rangle$, with a reward that might be negative (as well as non-negative). Similarly with orders: discarding an order for a product, which is labeled as item 
type $\ell$, is a matching $\langle \ell \rangle$ with the corresponding (most likely, negative) reward. 
Expediting component delivery can be included as well. Suppose matching $\langle 1,2,3,9 \rangle$
corresponds to product $9$ assembled from (one unit of) parts $1,2,3$, with the reward $20$.
However, the system has an option of expediting component $2$, and receive it immediately, at the 
cost of $15$. Then, assembling product $9$ from already available components $1$ and $3$, and expedited
component $2$, can be modeled as a matching $\langle 1,3,9 \rangle$ with reward $20-15=5$.
(Another, more natural, way to model expediting of item $2$ is to treat it as a ``matching,''
requiring $-1$ type-2 items, with the reward $-15$. See Remark~\ref{remark:negative-mu} above.)

This discussion illustrates the flexibility of our model {\em as long as the objective is to maximize average rewards
associated with actions}, such as matching, salvaging, expediting, etc. The model does {\em not} explicitly include holding costs. In Section~\ref{profit} we propose and discuss heuristic extensions of our scheme which do implicitly take holding costs into account.

\subsection{Virtual Matching System}
\label{secmatchvirtualdef}

We will propose a matching control scheme in Section \ref{secdynvirtual}, which 
in parallel to the physical system ``runs'' a virtual system, which determines the matching decisions for the
physical one.  The virtual matching system is defined as follows.

The virtual system has the same item types, set of matchings and arrival flows as the physical system. It is only different in that any matching can be activated at any time and the queues of the virtual system can be negative, as well as positive. Matchings in the virtual system are activated based on its own state, regardless of the state of physical system. The activated matchings in the virtual system become actual matchings in the physical system either immediately, or later in time, depending on availability of physical items. The virtual matchings, until they become actual ones, are called \emph{incomplete} matchings. Incomplete matchings wait in a queue, which lists the incomplete matchings (their identities $j$) in the order of arrival;
we denote the length of this queue by $\hat Q_{0}$. An incomplete matching becomes an actual one and leaves this queue when 
it is ``completed'' by all required physical items. 
(Incomplete matchings' queue, as we will see shortly, serves as the ``interface'' between the  virtual and physical systems. In our figures and plots it is shown as part of the physical system.)

%to become ``\emph{complete}'' with the future arrivals in FCFS order.

\subsection{Control of the Physical Matching System via Virtual System}
\label{secdynvirtual}

Denote by ${Q}(t)=({Q}_i(t),i\in\mathcal{I})$ and $\hat{Q}(t)=({Q}_i(t),i\in\mathcal{I})$
the vectors of queue lengths in the virtual and physical systems, respectively, at time $t$. 
In this paper we always assume that the system is initialized in a state such that all physical and virtual queues are 
zero, $Q_i(0)=\hat Q_i(0) = 0, ~\forall i \in \mathcal{I}$, and there are no incomplete matchings, $\hat Q_0(0)=0$.
%(We refer to this state as the {\em zero state.})
This means that the only feasible system states are those reachable from this ``zero-state.''
\iffalse
 we assume that the zero state is also reachable from any feasible state.  (If the process $Q(\cdot)\doteq (Q(t), ~t\ge 0)$ happens to be a Markov chain, the above assumptions simply mean that its state space consists of feasible states, and this Markov chain is irreducible.) 
 \fi

At time $t$ 
the following occurs sequentially:
\begin{enumerate}[label=(\roman*),font=\itshape]
\item A new matching is chosen in the virtual system based on $Q(t)$. (We will give a specific rule in Section~\ref{secmatchmodelalg}.) 
If it is a non-empty matching $j$, then the virtual queues are updated as $Q := Q -\mu(j)$, and a new type $j$ incomplete matching is created and placed at the end of the (incomplete matchings') queue; so that $\hat Q_0 := \hat Q_0+1$.
	\item The incomplete matchings' queue is scanned in FCFS order, to find the first incomplete matching $j'$, which can be completed, i.e. such that $\hat Q(t)\geq \mu(j')$. If such matching $j'$ is found,
it is completed, i.e. it is removed from the incomplete matchings' queue (so that $\hat Q_0 := \hat Q_0 - 1$), a physical matching $j'$ is created, and the corresponding number of physical items leaves the system, $\hat{Q} := \hat{Q} -\mu(j')$.
	\item Both $Q$ and $\hat{Q}$ are increased as: $Q :=Q + \lambda(t)$,
$\hat Q := \hat Q + \lambda(t)$; here $\lambda(t)=(\lambda_i(t),i\in\mathcal{I})$ is the random vector of arrivals of different types at $t$.
	
\end{enumerate}

According to steps \textit{(i)}-\textit{(iii)} above, if matching $j\in\mathcal{J}$ is chosen in the virtual system at time $t$, the virtual queues change as follows:
\begin{equation}
{Q}(t+1)={Q}(t)+\lambda(t)-\mu(j).\label{matchdynn2}
\end{equation}
The evolution of the physical queues, if matching $j'\in\mathcal{J}$ is completed is:
$$
\hat{Q}(t+1)=\hat{Q}(t)+\lambda(t)-\mu(j')
$$

\begin{figure}[h t!]
\centering
	\includegraphics[width=0.45\textwidth]{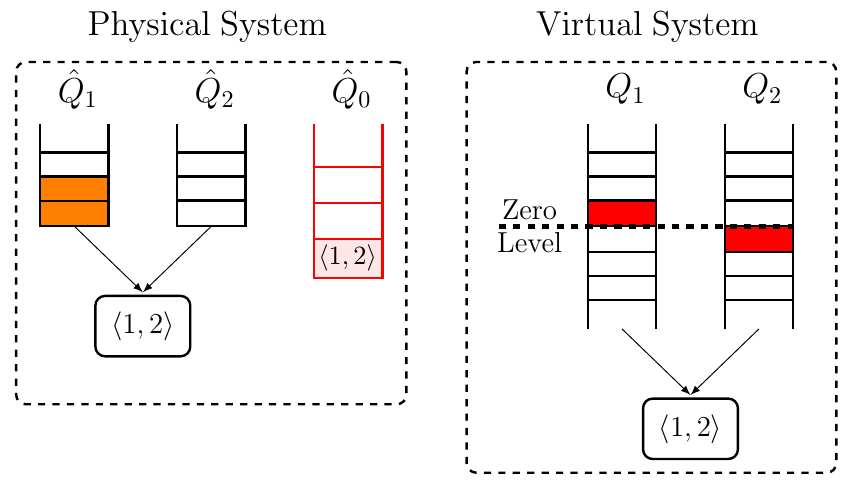}
\caption{An example of the physical and virtual matching systems}
\label{fig:virtual}
\end{figure}

Recall that we only consider feasible states of the queues -- those reachable from the state where all virtual and physical queues are zero. Then we can make the following 
observations for the control scheme described above. For illustration, we will use 
Figure~\ref{fig:virtual} showing a physical matching system with two item types and one possible matching and its corresponding virtual system. The system state shown on Figure~\ref{fig:virtual} is such that: (a) in the physical system there are two type 1 items and no type 2 items; (b) the queue lengths in the virtual system are $Q_1(t) = 1$, $Q_2(t) = -1$; (c) there is one incomplete matching 
$\langle1,2\rangle$, which is incomplete because, while there is a type 1 item in the physical system (to complete it), there is no available (physical) type 2 item. (Note that at this point we did not specified yet the matching rule(s) for the virtual system -- this will be done in Section~\ref{secmatchmodelalg}. So, the state on Figure~\ref{fig:virtual} only illustrates the relation between virtual and physical systems, not a specific matching rule.) 

In a general system, if ${Q}_i(t) <0$, then $Q^-_i(t)=|{Q}_i(t)|$ is the current shortage of type $i$ items for completing all incomplete matchings. (On Figure~\ref{fig:virtual}, ${Q}_2=-1$ indicates the shortage of one type 2 item for completion of the incomplete matching $\langle1,2\rangle$.)
If ${Q}_i(t) \ge 0$, then $Q^+_i(t)=Q_i(t)$ is the current surplus of type $i$ items,
beyond what is needed for completing all incomplete matchings.
(On Figure~\ref{fig:virtual}, ${Q}_1=1$ indicates that there is one type 1 item in addition 
to one type 1 item which can be used for completion of  the incomplete matching $\langle1,2\rangle$.)

In addition to notations $Q(t)$ and $\hat Q(t)$, let us denote by $\hat{\mathcal{Q}}_0(t)$ the state (list) of all incomplete matchings at time $t$.

The following simple Proposition~\ref{prop-queue-relation} gives a total queue length bound 
\eqn{eq-que-rel}
for the physical system in terms of the virtual one. This bound does not require any additional assumptions.
Statements (ii) and (iii) of the proposition involve the notion of stochastic stability, which means positive recurrence of a Markov process.
To keep the exposition simple, assume that the process 
%$Q(\cdot)\doteq (Q(t), t\ge 0)$,
$(Q(t), t\ge 0)$,
describing the evolution of the virtual system, and the process 
%$[Q(\cdot), \hat Q(\cdot), \hat Q_0(\cdot))] \doteq [(Q(t), t\ge 0), (\hat Q(t), t\ge 0), (\hat Q_0(t), t\ge 0)]$
$ [(Q(t), t\ge 0), (\hat Q(t), t\ge 0), (\hat{\mathcal{Q}}_0(t), t\ge 0)]$
describing the evolution of the entire system, are countable-state-space Markov chains.
(This is the case, for example, under the virtual system matching algorithm that we propose below in Section~\ref{secmatchmodelalg}, and under linear utility function $G$.)

\iffalse
Denote by ${Q}(t)=({Q}_i(t),i\in\mathcal{I})$ and $\hat{Q}(t)=({Q}_i(t),i\in\mathcal{I})$
the vectors of queue lengths in the virtual and physical systems, respectively, at time $t$. 
In this paper we always assume that the system is initialized in a state such that all physical and virtual queues are 
zero, $Q_i(0)=\hat Q_i(0) = 0, ~\forall i \in \mathcal{I}$, and there are no incomplete matchings, $\hat Q_0(0)=0$.
%(We refer to this state as the {\em zero state.})
This means that the only feasible system states are those reachable from this ``zero-state.''
\iffalse
 we assume that the zero state is also reachable from any feasible state.  (If the process $Q(\cdot)\doteq (Q(t), ~t\ge 0)$ happens to be a Markov chain, the above assumptions simply mean that its state space consists of feasible states, and this Markov chain is irreducible.) 
 \fi
 \fi

\begin{proposition}
\label{prop-queue-relation}
(i) At any $t\ge 0$, the following relation between physical and virtual queues holds:
\beql{eq-que-rel}
\hat Q_0(t) \le \sum_i Q_i^-(t), ~~~
\sum_i \hat Q_i(t) \le \sum_i Q_i^+(t) + \mu^* \sum_i Q_i^-(t) \le \mu^* \sum_i |Q_i(t)|,
\end{equation}
where $\mu^*\doteq \max_j \sum_i \mu_i(j)$.\\
(ii) Stochastic stability of $(Q(t), t\ge 0)$ implies that of $ [(Q(t), t\ge 0), (\hat Q(t), t\ge 0), (\hat{\mathcal{Q}}_0(t), t\ge 0)]$.\\
(iii) If $(Q(t), t\ge 0)$ is stochastically stable, then the steady-state average rates at which different matchings are activated are the same in the physical and virtual systems.
\end{proposition}

%\proof[Proof of something]
\proof
(i) Clearly, for all $i\in \mathcal{I}$ at all times, $Q_i(t) \le \hat Q_i(t)$. Note that the total shortage of items of all types for the completion of all incomplete matchings is $\sum_i Q_i^-(t)$; this means, in particular, 
that the total number of incomplete matchings is upper bounded as $\hat Q_0(t) \le \sum_i Q_i^-(t)$. The total number
of physical items in the system, $\sum_i \hat Q_i(t)$, can be partitioned into those that are ready to be used for
completion of incomplete matchings and the ``surplus'' items; the number of the former is upper bounded
by $\mu^* \hat Q_0(t)$;
the number of the latter is equal to $\sum_i Q_i^+(t)$. This implies the second part of \eqn{eq-que-rel}.

(ii) Follows from (i).

(iii) Follows from (ii).
\endproof

\begin{remark}
\label{remark:gen-m} \normalfont
If $m \ge 1$ matchings can be done after each arrival, the sequence of steps (i)-(iii) above is repeated $m$ times.
\end{remark}

\subsection{Asymptotically Optimal Matching Algorithm for The Virtual System}
\label{secmatchmodelalg}

We now specify the algorithm to be used for the control of the virtual system. This algorithms will be proved to be asymptotically optimal for the virtual system, and then (by Proposition~\ref{prop-queue-relation}) for the physical system as well -- see Remark~\ref{rem-both-optimal} below.

\begin{algorithm}[h]
\caption{Matching Algorithm for the Virtual System}
\label{alg-virtual}
\normalsize 
Let a (small) parameter $\beta>0$ be fixed.
At each time $t=1,2,\cdots$, activate matching
\begin{equation}
 j(t)\in\argmax_{j\in \mathcal{J}} \left[ (\partial G(X(t))/\partial X_j)~w_j+\sum_{i\in{\mathcal{I}}}\beta {Q}_i(t)\: \mu_i(j)\right], \label{rule3}
\end{equation}
where running average values $X_j(t)$ (of the rewards obtained by activation of different matchings $j$) are updated as follows:
\begin{eqnarray}
X_{j(t)}(t+1)&=&(1-\beta) X_{j(t)}(t)+\beta\, w_{j(t)},\\
 X_{j}(t+1)&=&(1-\beta) X_{j}(t),\qquad j\neq j(t),
\end{eqnarray}
and ${Q}_i(t)$ is updated according to rule (\ref{matchdynn2}) for all $i\in\mathcal{I}$.
\end{algorithm}

Note that if the function $G$ is linear, say $G(X) = \sum_j X_j$, then the partial derivatives in \eqn{rule3}
are constant, and rule \eqn{rule3} becomes simply
\beql{rule3-linear}
j(t)\in\argmax_{j\in \mathcal{J}} \left[ w_j+\sum_{i\in{\mathcal{I}}}\beta {Q}_i(t)\: \mu_i(j)\right].
\end{equation}
Moreover, in this case the algorithm does {\em not} need to keep track of the averages $X_j(t)$. As a result,
both processes $(Q(t), t\ge 0)$ and
$ [(Q(t), t\ge 0), (\hat Q(t), t\ge 0), (\hat{\mathcal{Q}}_0(t), t\ge 0)]$
are countable-state-space Markov chains.

Consider the following Assumption~\ref{assumption-matching} on the model structure. It is stated informally -- its precise meaning will be given (in a more general context) later in Assumption~\ref{assumption 1} (Section~\ref{sec3network}).
Also, in Section~\ref{sec-assump-nonrestrict} we explain why this assumption is non-restrictive.

\vspace{.2cm}
\begin{assumption}\label{assumption-matching}
For any subset $\bar{\mathcal{I}}\subseteq \mathcal{I}$, there exists 
 a matching activation strategy, under which the long-term average drift of queues $i\in\bar{\mathcal{I}}$
is strictly positive and the long-term average drift of queues $i \not\in \bar{\mathcal{I}}$
is strictly negative.
\end{assumption}

When parameter $\beta$ is small, then the {\em running average $X_j(t)$
is (one notion of) a long-term average rate} at which rewards due to matching $j$ are generated.
(See Section~\ref{subsec-egpd}.)
We will prove in Section~\ref{sec3network} (as a corollary of Theorem~\ref{zTheorem 2}) 
that, {\em under Assumption~\ref{assumption-matching}, Algorithm~\ref{alg-virtual} is asymptotically optimal} 
in the following sense. (It is described here informally -- the formal result is Theorem~\ref{zTheorem 2}, 
for the more general model in Section~\ref{sec3network}.)
Let $V$ be the set of those long-term rate vectors $X$  that are achievable (by some control strategy) subject to the stability of the queues, and let $V^*$ be its optimal subset, $V^* = \argmax_{X\in V} G(X)$.
%Let $V^*$ be set of those long-term rate vectors $X$ that achieve the maximum possible $G(X)$ (over all control strategies), subject to stability of the queues;
Then, when $\beta$ is small, $X(t) \to V^*$
as $t\to\infty$.

Suppose now that the system process is Markov under Algorithm~\ref{alg-virtual}
(as is the case when function $G$ is linear).
Then  Assumption~\ref{assumption-matching}
ensures the process stability (for example, by the argument described in Section 4.9 in \cite{stolyar2005maximizing}). In this case the steady-state average rewards
(due to different matchings) $u=(u_1,\ldots,u_J)$ are well defined.
If the process is in stationary regime, then obviously $\mathbb{E} X(t)=u$.
Furthermore, the asymptotic optimality of Algorithm~\ref{alg-virtual} in the sense described above,
can be used to show that, as  $\beta\to 0$,
the vector $u$ converges to the optimal set $V^*$ (see Section 4.9 in \cite{stolyar2005maximizing}).

\begin{remark} \normalfont
\label{rem-both-optimal}
If Algorithm~\ref{alg-virtual} is asymptotically optimal for the virtual system, then under our scheme it is also asymptotically optimal for the physical system. Indeed, 
the physical and virtual systems have the same set $V$ of achievable long-term rate vectors $X$ (subject  to the stability of the queues). This is because any $X$ achievable in the virtual system is achievable in the physical system as well
(by our scheme, for which we have Proposition~\ref{prop-queue-relation}), and vice versa because obviously any control of the physical system can be applied to the virtual system. Therefore, under our scheme, if the virtual system produces (in the asymptotic limit) the optimal long-term rates $X \in V^*$, the same optimal rates are produced (by
Proposition~\ref{prop-queue-relation}) in the physical system.
\end{remark}

\subsection{Discussion of Algorithm~\ref{alg-virtual}}
\label{sec-virtual-alg-discussion}

\subsubsection{Basic intuition}

The key feature of the virtual system is that it has an option of creating matchings ``in advance,''
before all required physical items have arrived. These ``advance'' matchings are the ones we called incomplete.
Virtual queues keep track of the items' availability: recall that if $Q_i < 0$, $|Q_i|$ is the shortage of type $i$ items, and  if $Q_i \ge 0$, it is the surplus of type $i$ items. 

The intuition behind Algorithm~\ref{alg-virtual} is the same as for the GPD algorithm
in \cite{stolyar2005maximizing}
(and other related works -- see, e.g., \cite{St2005gpdgen} and references therein),
but our model is more general in that the queues may have any sign. 
%(For other related work, see, e.g., \cite{St2005gpdgen} and references therein.)
For simplicity of discussion, suppose the objective function is linear, $G(X)=\sum_j X_j$,
in which case Algorithm~\ref{alg-virtual} specializes to \eqn{rule3-linear}.
The rule ``tries'' to choose a matching $j$ which brings large reward $w_j$,
but at the same time it ``encourages'' the drift of the queues towards $0$.
%tries to minimize the increment of $\sum_i Q_i^2$.
Indeed, recall that activation of any matching can only decrease the virtual queues.
This means that the rule ``encourages'' the use of matchings that decrease positive $Q_i$'s as much as possible and decrease negative $Q_i$'s as little as possible; in other words, 
the rule encourages matchings requiring items of which there is a large surplus,
and discourages matchings requiring items of which there is already a large shortage -- this guarantees 
stability of the queues. When parameter $\beta$ is small, the virtual queues ``stabilize around correct levels'' --
positive or negative -- which allows rule \eqn{rule3-linear} to make ``correct'' decisions maximizing the average rewards.

\subsubsection{Assumption~\ref{assumption-matching} is non-restrictive}
\label{sec-assump-nonrestrict}

We now describe two common cases, in which Assumption~\ref{assumption-matching} holds. These two cases cover a very large number of applications.

{\em Case 1.}  Assumption~\ref{assumption-matching} holds automatically in the special case when 
for each item type $i$ there exists at least one matching requiring only type $i$ items (namely, with $\mu_i \ge 1$ and $\mu_\ell =0$ for $\ell \ne i$). In this case it suffices to pick any parameter $m$ (the number of matchings per each batch arrival) which is greater than $\mu^* \doteq \max_j \sum_i \mu_i(j)$. 
This special case is very common for the following reason, which we illustrate using 
the simple model in Figure~\ref{fig:virtual}.
If matching $\langle 1,2 \rangle$ is the only possible (besides the empty matching), the system is
unbalanced when the arrival rates are unequal, $\alpha_1 \ne \alpha_2$, and cannot be stable. 
(If items arrive one-by-one, this particular system obviously cannot be stable even if $\alpha_1 = \alpha_2$.
More generally, any system with one-by-one arrivals cannot be stable if its ``matching graph'' is bi-partite,
see \cite{mairesse2014stability}.)
This shows that many practical systems typically need the option of using ``single'' matchings $\langle i \rangle$ anyway (salvaging or discarding individual items), to ensure stability, and then 
Assumption~\ref{assumption-matching} holds. 

{\em Case 2.}  This case is more subtle. Suppose a system can potentially be made stable without requiring single-type matchings.
For example, consider the system in Figure~\ref{fig:virtual} in which the arrivals occur only in pairs $(1,2)$.
Suppose also that up to two matchings can be done upon each arrival ($m=2$).
On the face of it, Assumption~\ref{assumption-matching} does not hold for this system.
Indeed, the linear relation $Q_1(t)=Q_2(t)$ holds at all times and, therefore, it is impossible for $Q_1$ and $Q_2$ to have different average drifts, which is required under Assumption~\ref{assumption-matching}. However, consider the orthogonal change of coordinates, $\tilde Q_1 = Q_1+Q_2, \tilde Q_2 = Q_1 - Q_2$, with $\lambda(\cdot)$ and 
$\mu(\cdot)$ transformed accordingly. 
Then, $\tilde Q_2(t) \equiv 0$,  and the system can be considered as having only one queue $\tilde Q_1$. For the latter system Assumption~\ref{assumption-matching} does hold. Note that {\em the algorithm itself does not need to perform any change of coordinates -- it remains as is.} This situation is generic: if there is an inherent linear dependence between 
the queues, Assumption~\ref{assumption-matching} often holds for the system after an appropriate orthogonal change of coordinates. This is, in fact, the case for many bi-partite matching systems (with items arriving in pairs), including the one we consider later in Section~\ref{subsec-bipart-sim}.

To summarize the discussion in this subsection, Assumption~\ref{assumption-matching} is essentially the assumption that the system can be made stable, plus a very common condition that the queues ``can be moved in any direction'' {\em within the subspace of feasible queue states.}

%%%%%%%%%%%%%%%%%%%%%%%%%%%%%%%%%%%%%%%%%%%%%%%%%%%%%%%%%%%%%
\section{A General Network Model and EGPD Algorithm}\label{sec3network}

In this section we introduce the \emph{Extended Greedy Primal-Dual} (EGPD) algorithm for a general network model, which includes the matching system as a special case. This algorithm is a generalization of the GPD algorithm of \cite{stolyar2005maximizing} in the sense that queues at some network nodes, we call them \emph{free} nodes, are allowed to have any sign; as they evolve, 
these queues are ``free'' to change from positive to negative and vice versa.
The model in \cite{stolyar2005maximizing} is such that queues at all nodes are constrained to be non-negative -- in 
our model we call such nodes {\em constrained}.
 First, we will formally define the model and the underlying optimization problem in Sections \ref{secThemodel}-\ref{sec:optprob}. The optimization problem determines the best possible (under any control
algorithm) long-term drifts of the queues, which maximize the network ``utility'' subject to the 
condition that queue-drifts are zero at free nodes and are non-positive at constrained nodes;
the optimal solutions to this problem give the maximum possible network utility that can be achieved by
any network control strategy subject to stability of the queues.
We define the EGPD algorithm in Section \ref{secGPDalgorithm}. In Section \ref{secFSP}, we show that, as the algorithm parameter $\beta$ goes to $0$,  the ``fluid scaled'' version of the process converges to a random process with sample paths being what we define as EGPD-trajectories. In Section \ref{sec:globalattraction}
we prove asymptotic optimality of the EGPD-algorithm,
in the sense that EGPD-trajectories converge to the optimal set of the underlying optimization problem while
keeping all queues uniformly bounded;
in other words, EGPD-algorithm maximizes the system utility subject to stability.
Finally, in Section~\ref{secmatchmodelenhance} we show that Algorithm~\ref{alg-virtual} 
(Section~\ref{secmatchmodelalg})
for the virtual system of 
Section~\ref{secmatchvirtualdef} is a special case of EGPD.

A reader interested only in the application of EGPD algorithm to the dynamic matching model 
of Section~\ref{secmatchmodel} may wish to skip 
at first reading the proofs in Sections \ref{secFSP}-\ref{sec:globalattraction}.

\subsection{The Model}\label{secThemodel}
Consider a network consisting of a finite set of nodes $\mathcal{N}=\{1,2,\cdots,N\}$, $N\geq 1$. The nodes are of two different types: $N_1$ \emph{constrained} nodes form the set $\mathcal{N}^c=\{1,2,\cdots,N_1\}$ and $N_2=N-N_1$ \emph{free} nodes form $\mathcal{N}^f=\{N_1+1,N_1+2,\cdots,N\}$. Either $\mathcal{N}^c$ or $\mathcal{N}^f$ is allowed to be an empty set. There is a queue associated with each node, where we denote by $Q_n(t)$ the queue length of node $n\in\mathcal{N}$ at time $t$ and we will denote $Q(t) = (Q_n(t), n\in\mathcal{N})$. The queue length of node $n\in\mathcal{N}^c$ is always \emph{non-negative}, but node $n\in\mathcal{N}^f$ can have queue length of any sign.

The system operates in discrete time $t=1,2,\cdots$. (By convention, we identify an integer time $t$ with unit time interval [t,t+1), which is usually referred to as time slot $t$.) A finite number of \emph{controls} is available, where we denote by $K$ the set of controls. With activation of control $k\in K$ at time $t$, the following occurs sequentially:
\begin{enumerate}[label=(\roman*),font=\itshape]

\item A certain (non-random) real amount (``number'') $\mu_n(k)\geq 0$ of items is removed from queue $n$ and leaves the network. Queues in constrained nodes cannot go below zero; 
so if $Q_n(t)\leq \mu_n(k)$, the entire content of queue $n$ is removed.

\item A random (bounded) real amount (``number'') $\lambda_n(k,t) \ge 0$ of items  enters each node $n\in\mathcal{N}$, where 
$\lambda(k,t) = (\lambda_n(k,t), n\in\mathcal{N})$ are i.i.d. in time, with generic random variable denoted 
$\lambda(k) = (\lambda_n(k), n\in\mathcal{N})$.

\end{enumerate}
According to steps \textit{(i)} and \textit{(ii)}, the queue update rules for constrained and free nodes,
given control $k$ is chosen at time $t$, 
 are as follows:
\begin{eqnarray}
&&Q_n(t+1)=[Q_n(t) - \mu_n(k)] \vee 0 +  \lambda_n(k,t) ,\quad n\in\mathcal{N}^c\label{modupd1}\\
&&Q_n(t+1)=Q_n(t) - \mu_n(k) + \lambda_n(k,t) ,\quad n\in\mathcal{N}^f \label{modupd2}.
\end{eqnarray}

\subsection{System Rate Region}\label{rateregion}
For each $k\in K$ and time $t$, consider the random vector $b(k,t)=(b_n(k,t), n\in\mathcal{N})$ equal in distribution to $\lambda(k) - \mu(k)$. Clearly, $b(k,t)$ is equal to the random vector of queue increments $Q(t+1)-Q(t)$ provided that control $k$ is chosen at time $t$ and assuming $Q_n(t)\geq \mu_n(k)$ for all $n\in\mathcal{N}^c$. We call
components of $b(k,t)$ the \emph{nominal increments} of queues upon control $k$ at time $t$. Let $k(t)$ denote the control chosen at time $t$ by a given control policy. 

Informally speaking, the finite-dimensional convex compact rate region $V\subset \mathbb{R}^{N}$ is defined as the set of all possible long-term average values of $b(k(t),t)$, which can be induced by different control policies.
Formal definition of the rate region is as follows. 

For each $k\in K$, denote by $\overline{b}(k)=\mathbb Eb(k,t) $ the drift of queue lengths upon control $k$ (at any time $t$ when control $k$ is activated). For a fixed probability distribution $\phi=(\phi_k, k\in K)$ (with $\phi_k\geq 0$ and $\sum_{k\in K}{\phi_k}=1$)  consider the vector
\begin{equation}
v(\phi)= \sum_{k\in K}{\phi_k \overline{b}(k)}.
\end{equation}
If we interpret $\phi_k$ as the long-term average fraction of time slots when control $k$ is chosen from the set of controls $K$, then $v(\phi)$ corresponds to the vector of long-term average drifts of $Q(t)$, assuming that the queues in the constrained nodes never hit zero. Then the system rate region $V$ is defined as the set of all possible vectors $v(\phi)$ corresponding to all possible $\phi$.

\subsection{Underlying Optimization Problem}\label{sec:optprob}

Consider an open convex set $\tilde{V}\subseteq \mathbb{R}^{N}$ such that $\tilde{V}\supseteq V$. Consider a concave continuously differentiable \emph{utility} function $H:\tilde{V}\rightarrow \mathbb R$ and the following optimization problem:
\begin{eqnarray}
\max_{v\in V} &&H(v)\label{qa1}\\
\textit{s.t. } &&v_n\in \mathbb{R}^{}_-,\:\forall\: n\in\mathcal{N}^c\nonumber\\
 &&v_n=0,\quad\forall\: n\in\mathcal{N}^f. \nonumber
\end{eqnarray}

\begin{assumption}\label{qnonemptyv}
Optimization problem {\normalfont(\ref{qa1})} is feasible, i.e.
\begin{equation}
\{v\in V: v_n\in\mathbb{R}^{}_- ,\forall n\in \mathcal{N}^c \textit{ and }v_n=0 ,\forall n\in \mathcal{N}^f \}\neq \varnothing.
\end{equation}
\end{assumption}

If Assumption~\ref{qnonemptyv} holds, we denote by $V^*\subseteq V$ the set of optimal solutions
of \eqn{qa1}. The dual to optimization problem (\ref{qa1}) is
\begin{equation}
\min_{(y_n\in \mathbb{R}^{}_+,n\in\mathcal{N}^c),(y_n\in \mathbb{R}^{},n\in\mathcal{N}^f)}\left(\max_{v\in V}\left(H(v)-y\cdot v\right)\right),\label{qa2}
\end{equation}
and we denote by $Q^*$ the closed convex set of optimal solutions $q^*\in \mathbb R^{N_1}_+\times \mathbb R^{N_2}$ of problem (\ref{qa2}). For any $v^*\in V^*$ and any $q^*\in Q^*$, the compementary slackness condition holds:
\beql{eq-comp-slack}
q^* \cdot v^* = 0.
 \end{equation}
%\vspace{.2cm}

In Section \ref{secGPDalgorithm}, we will introduce an algorithm, which is asymptotically optimal under the following assumption, which is stronger than Assumption~\ref{qnonemptyv}.

\vspace{.3cm}
\begin{assumption}\label{assumption 1}
For any subset $\bar{\mathcal{N}}^f\subseteq \mathcal{N}^f$, there exists $v\in V$ such that $v_n>0$ for $n\in\bar{\mathcal{N}}^f$ and $v_n<0$ for $n\not\in\bar{\mathcal{N}}^f$.
\end{assumption}
\noindent Assumption \ref{assumption 1} means that there always exists a control policy which provides, simultaneously, a \textit{strictly negative average drift} to all the constrained node queues and \textit{non-zero average drifts} toward zero for all free node queues.

Note that under Assumption \ref{assumption 1}, the set $Q^*$ is compact. Indeed, the optimal value of the problem (\ref{qa1}) is equal to
\begin{equation}
H(v^*)= \max_{v\in V} \left(H(v)-q^*\cdot v\right)\label{qaa2}
\end{equation}
for any $v^*\in V^*$ and any $q^*\in Q^*$. Set $Q^*$ must be bounded, because otherwise, from Assumption \ref{assumption 1}, there would exist $v\in V$ such that $v_n<0$ for all nodes with $q_n\geq0$, and $v_n>0$ for all nodes with $q_n<0$. Then we can arbitrarily increase the RHS of (\ref{qaa2}) by choosing $q^*\in Q^*$ with large $|q_n^*|$.

The problem that we are going to address is as follows. Let $X$ denote a long-term average value of $b(k(t),t)$ under a given dynamic control policy, that is, a policy of choosing $k(t)$ depending on the system state. We are interested in finding a dynamic control policy such that when optimization problem (\ref{qa1}) is feasible, and moreover, the stronger Assumption~\ref{assumption 1} holds, the corresponding $X$ is close to $V^*$, while the system queues remain stochastically stable. 

\subsection{Extended Greedy Primal-Dual Algorithm}\label{secGPDalgorithm}
\label{subsec-egpd}

Consider the following control policy: 
\begin{algorithm}[H]

\caption{EGPD algorithm for the general network model}
\label{alg-egpd}
%\begin{algorithmic}[1]
\normalsize 
At time $t=1,2,\cdots$, choose a control
\begin{equation}
k(t)\in\argmax_{k \in K} \left[\nabla H(X(t))-\beta Q(t) \right]\cdot\overline{b}(k),\label{eq:matchargmax}
\end{equation}
where $\beta>0$ is a small parameter. Here $X(t)$ is the running average of $b(k(t),t)$, updated as follows:
\begin{equation}
X(t+1)=(1-\beta) X(t)+\beta\,b(k(t),t)\label{eq:matchupdateruleforx}
\end{equation}
and $Q(t)$ is updated according to (\ref{modupd1}) and (\ref{modupd2}).
\end{algorithm}
\vspace{-.3cm}
The initial condition is $X(0)\in \tilde{V}$. Note that such initial condition and update rule (\ref{eq:matchupdateruleforx}) imply that $X(t)\in\tilde{V}$ for all $t\geq 0$. Hence the system evolution is well-defined for all $t\geq 0$, since the gradient and argmax in (\ref{eq:matchargmax}) are well-defined.

Also note that, if $0<\beta<1$, then for $t\ge 1$
$$
X(t) = \sum_{\tau=0}^{t-1} \beta (1-\beta)^{t-1-\tau} b(k(\tau),\tau) +
(1-\beta)^t X(0).
$$
Therefore, when $t$ is large, $X(t)$ is essentially the 
geometric average of  values of $b(k(\tau),\tau)$ up to time $t-1$.
When $t$ is large and $\beta>0$ is small, $X(t)$ is (one notion of) the long-term average of  values of $b(k(\tau),\tau)$ 
up to time $t-1$.

\subsection{Asymptotic Regime and Fluid Limit}\label{secFSP}
We define \textit{EGPD-trajectory} as a pair of absolutely continuous functions $(x,q) = ((x(t), t \geq 0), (q(t), t \geq 0))$, each taking values in $\mathbb R^{N}$ and satisfying the following conditions:

\textit{(i)} For all $t\geq 0$,
\begin{equation}
x(t)\in\tilde{V}\label{eq:FSP-19-1}
\end{equation}
and for almost all $t\geq 0$,
\begin{equation}
x'(t)=v(t)-x(t),\label{eq:FSP-19-2}
\end{equation}
where 
\begin{equation}
v(t)\in\argmax_{v\in V} [\nabla H(x(t))-q(t)]\cdot v\label{eq:FSP-19-4}.
\end{equation}
\textit{(ii)} We have
\begin{eqnarray}
%&q_n(0)\geq 0,\, n\in \mathcal{N}^c\label{eq:FSP-19-1-5}\\
&q_n(t)\geq 0,\, \forall t\geq 0 ,\; n\in\mathcal{N}^c\label{qab1}\\
&q'_n(t)=[v_n(t)]^+_{q_n(t)},\;\textit{ a.e. in }t\geq 0,\, n\in\mathcal{N}^c\label{qab2}\\
&q'_n(t)=v_n(t),\;\textit{ a.e. in }t\geq 0,\, n\in\mathcal{N}^f\label{qab3}
\end{eqnarray}
Functions $x(t)$ and $q(t)$ are dynamically changing primal and dual variables, respectively, for problems (\ref{qa1}) and (\ref{qa2}), which arise as asymptotic limits of the fluid scaled version of the process as described next.

Consider a sequence of processes $(X^\beta,Q^\beta)$, indexed by a parameter $\beta$, where $\beta\downarrow 0$ along a sequence $\mathcal{B}=\{\beta_j\}_{j=1}^\infty$ with $\beta_j>0$ for all $j$. The initial state $(X^\beta(0),Q^\beta(0))\in\tilde{V}$ is fixed for each $\beta\in\mathcal{B}$. (The processes and variables associated with a fixed parameter $\beta$ will be supplied by superscript $\beta$.)

We need to augment the definition of the process. Let us assume $X^\beta(t)$ and $Q^\beta(t)$ are functions defined on $t\in\mathbb{R_+}$ and constant within each time slot $[l,l+1)$, $l=0,1,2,\cdots$. Thus for each $\beta$, consider the (continuous-time) process $Z^\beta=(X^\beta,Q^\beta)$, where
\begin{eqnarray}
X^\beta=(X^\beta(t)=(X_n^\beta(t),n\in\mathcal{N}), t\geq 0),\\
Q^\beta=(Q^\beta(t)=(Q_n^\beta(t),n\in\mathcal{N}), t\geq 0).
\end{eqnarray}
For each $\beta$, 
\begin{equation}
z^\beta= (x^\beta,q^\beta)
\end{equation}
is the fluid scaled version of process $Z^\beta$, obtained by
\begin{equation}
x^\beta=X^\beta(t/\beta),\quad q^\beta=\beta Q^\beta(t/\beta).\label{eq:rescale}
\end{equation}
The following theorem is straightforward modification of Theorem 3 in \cite{stolyar2005maximizing}, which we present without proof.

\vspace{.3cm}
\begin{theorem}\label{thm3}
Consider a sequence of process $\{z^\beta\}$ with $\beta\downarrow 0$ along set $\mathcal{B}$. Each process is considered as a random element in the Skorohod space of RCLL (``right continuous with left limits'') functions. Assume that $z^\beta(0)\rightarrow z(0)$, where $z(0)$ is a fixed vector in $\mathbb R^{2N}$ such that $X(0)\in \tilde{V}$. Then, the sequence $\{z^\beta\}$ is relatively compact and any weak limit of this sequence (i.e a process obtained as the weak limit of a subsequence of $\{z^\beta\}$) is a process with sample paths $z$ being EGPD-trajectories (with initial state $z(0)$) with probability 1.
\end{theorem}

\subsection{Global Attraction Result}\label{sec:globalattraction}
The following theorem is the main result of this section which shows the convergence of EGPD-trajectories to the saddle set $V^*\times Q^*$.
\vspace{.3cm}
\begin{theorem}\label{zTheorem 2}
Under Assumption \ref{assumption 1}, the following holds:
\begin{enumerate}[label=(\roman*),font=\itshape]
\item For any EGPD-trajectory $(x,q)$, as $t\rightarrow \infty$,
\begin{eqnarray}
	&x(t)\rightarrow V^*,\label{eq:converge}\\
	&q(t)\rightarrow q^{*}\textit{, for some } q^{*}\in Q^*\label{eq:qconverge}.
	\end{eqnarray}
\item Let some compact subsets $V^{\Box}\subset\tilde{V}$ and $Q^{\Box}\subset \mathbb{R}^{N_1}_+\times \mathbb{R}^{N_2} $ be fixed. Then, the convergence 
\begin{equation}
(x(t),q(t))\rightarrow V^*\times Q^*, ~~ t\rightarrow \infty,
\end{equation}
is uniform across all EGPD-trajectories with initial states $(x(0),q(0))\in V^{\Box}\times Q^{\Box}$.
\end{enumerate}
\end{theorem}

The proof of Theorem~\ref{zTheorem 2} is a generalization of that of Theorem 2 in \cite{stolyar2005maximizing} --
all steps of the latter are extended to our more general setting. For this reason
we will not give a complete proof of Theorem~\ref{zTheorem 2} in this paper, because it is lengthy. 
Instead, we demonstrate the key points involved in the generalization, by proving in this section the convergence
\eqn{eq:converge}
%Theorem~\ref{zTheorem 2} 
for the special case when $x(0)\in V$ and $H(\cdot)$ is \emph{strictly} concave.

Consider a fixed EGPD-trajectory $(x,q)$. The property 
\begin{eqnarray}
\rho(x(t),V)\leq \rho(x(0),V) e^{-t}, \: t\geq 0
\label{rho1}
\end{eqnarray}
holds regardless of Assumptions \ref{qnonemptyv} or \ref{assumption 1} (cf. Lemma 20 in \cite{stolyar2005maximizing}). This shows that entire trajectory $(x(t),t\geq 0)$ is contained within $V$. This fact implies that 
$\sup_{t\ge 0} \|\nabla H(x(t))\| < \infty$.
%$\nabla H(x(t))$ is uniformly bounded for all $t\geq 0$. 

A time point $t\geq 0$ is called ``\textit{regular}'' if conditions (\ref{eq:FSP-19-1})-(\ref{eq:FSP-19-4}) are satisfied and proper derivatives $x'(t)$, $q'(t)$ and $f'(t)$ exist. Almost all $t$ are regular.

Let us introduce the following function:
\begin{eqnarray}
F(v,y)=H(v)-\frac{1}{2}\sum_{n\in \mathcal{N}} y_n^2,\quad v\in \tilde{V},\nonumber\label{eq23}
&\:y_n\in \mathbb{R}^{}_+\text{ for } n\in\mathcal{N}^c, \nonumber
&\:y_n\in \mathbb{R}^{}\text{ for }n\in\mathcal{N}^f.\nonumber
\end{eqnarray}

\begin{lemma}
\label{lemma3}
The trajectory $(q(t), t\geq 0)$ is such that 
\begin{equation}
\sup_{t\geq 0} \left\|q(t)\right\|<\infty.
\end{equation}
\end{lemma}

\proof Within this proof, we say that a vector-function (or scalar function) $\alpha(t), ~t\ge 0,$ is {\em uniformly bounded} if $\sup_{t\ge 0} \|\alpha(t)\| < \infty$. 
By Assumption \ref{assumption 1}, the following holds for some fixed number $\delta>0$. For any  $t\ge 0$, there exists $\xi =(\xi_n,\, n\in \mathcal{N})\in V$ such that for any $n$, $|\xi_n|\ge \delta$, $\xi_n>0$ if $q_n<0$, and $\xi_n<0$ if $q_n\geq 0$. Then for any regular $t \ge 0$ (and a corresponding $\xi$) 
 we have:
\begin{eqnarray}
\frac{d}{dt}F(x(t),q(t))&=& [\nabla H(x(t))-q(t)]\cdot v(t)-\nabla H(x(t))\cdot x(t)\nonumber\\
&\geq& [\nabla H(x(t))-q(t)]\cdot \xi-\nabla H(x(t))\cdot x(t)\nonumber\\
&=& -\sum_{n\in\mathcal{N}}\xi_n q_n(t)+\nabla H(x(t))\cdot(\xi-x(t))\nonumber\\
&\geq&\delta \sum_{n\in\mathcal{N}} |q_n(t)|+\nabla H(x(t))\cdot(\xi-x(t))\label{q32}
\end{eqnarray}
Since $\nabla H(x(t))$ and $x(t)$ are uniformly bounded, so it the second term in (\ref{q32}).
When $\|q(t)\|$ is large, the first term in (\ref{q32}) is large positive. 
We conclude that  $({d}\slash{dt}) F(x(t),q(t)) \geq\epsilon_1>0$ as long as $\left\|q(t)\right\|\geq C_1>0$, for some fixed constants $\epsilon_1$ and $C_1$. Since $H(x(t))$ is uniformly bounded, we can pick $C_2>0$ sufficiently large so that
$F(x(t),q(t))\leq -C_2$ implies $\left\|q(t)\right\|\geq C_1$ and then $({d}\slash{dt}) F(x(t),q(t)) \geq\epsilon_1>0$.
We then conclude that $\liminf_{t\to\infty} F(x(t),q(t)) \ge -C_2$ and, therefore, $\inf_{t\ge 0} F(x(t),q(t)) > -\infty$.
The latter (along with the uniform boundedness of  $H(x(t))$) implies that  $q(t)$ is uniformly bounded.
\endproof

\begin{lemma}\label{Lemma4}
\vspace{.3cm}
For any EGPD-trajectory, at any regular time $t\geq 0$,
\begin{equation}
\frac{d}{dt}F(x(t),q(t))=\nabla H(x(t))\cdot(v(t)-x(t))-q(t)\cdot v(t)\label{qa31}
\end{equation}
and
\begin{equation}
v(t)\in\argmax\limits_{v\in V} \nabla H(x(t))\cdot(v-x(t))-q(t)\cdot v\label{qa32}
\end{equation} 
Furthermore, if Assumption {\normalfont\ref{qnonemptyv}} holds,
\begin{equation}
\frac{d}{dt}F(x(t),q(t))\geq\nabla H(x(t))\cdot(v^*-x(t))\geq H(v^*)-H(x(t)).\label{qa33}
\end{equation}
\end{lemma}
\proof
Noting $q'_n(t)=v_n(t)$ and $v_n^* = 0$, for any $n\in\mathcal{N}^f$, every step of the proof is analogous to that of Lemma 3 in \cite{stolyar2005maximizing}.
\endproof

Select an arbitrary point $q^{*}\in Q^*$ and associate it with the following function
\begin{eqnarray}
F^*(v,y)=H^*(v)-\frac{1}{2}\sum_{n\in \mathcal{N}}{(y_n-q_n^{*})^2},\quad v\in\tilde{V}, \nonumber
&\:y_n\in \mathbb{R}^{}_+\text{ for }n\in\mathcal{N}^c, \nonumber
&\:y_n\in \mathbb{R}^{}\text{ for }n\in\mathcal{N}^f,\label{eq27}
\end{eqnarray}
where 
\begin{equation}
H^*(v)=H(v)-q^{*}\cdot v\label{b01}\nonumber
\end{equation}
is the Lagrangian of problem (\ref{qa1}) with the dual variable equal to $q^{*}\in Q^*$. 
Having strictly concave $H(\cdot)$ implies that $H^*(\cdot)$ is also a strictly concave function and
\begin{equation}
v^*=\argmax_{v\in V} H^*(v)
\end{equation} 
is the unique optimal solution. Note that $\nabla H^*(v) = \nabla H(v)-q^*$. 
\vspace{.3cm}
\begin{lemma}\label{lemma5}
Consider $F^*(\cdot,\cdot)$ associated with an arbitrary $q^{*}\in Q^*$. Then for all (regular) $t\geq 0$,
\begin{equation}
\frac{d}{dt}F^*(x(t),q(t))\geq\left[\nabla H(x(t))-q^{*}\right]\cdot(v(t)-x(t))- (q(t)-q^{*})\cdot v(t)\label{a35}
\end{equation}
and 
\begin{equation}
x(t)\in V \emph{ implies }\frac{d}{dt}F^*(x(t),q(t))\geq 0.\label{a36}
\end{equation}
\end{lemma}
\proof
The proof is analogous to that of Lemma 5 in \cite{stolyar2005maximizing}. The only difference is the existence of free nodes, where we can easily validate this Lemma by using $q'_n(t)=v_n(t)$ and $v_n^*=0$ for any $n\in\mathcal{N}^f$.
\endproof

\proof[Proof of Theorem \ref{zTheorem 2}]
The convergence \eqn{eq:converge} follows from an inequality that we first derive. For any (regular) $t\geq 0$,
\begin{eqnarray}
&\frac{d}{dt}F^*(x(t),q(t))\geq  (\nabla H(x(t))-q^*)\cdot(v(t)-x(t))-(q(t)-q^*)\cdot v(t)\nonumber\\
&=  (\nabla H(x(t))-q^*)\cdot(v^*-x(t))-(q(t)-q^*)\cdot v(t) + (\nabla H(x(t))-q^*)\cdot(v(t)-v^*)\nonumber\\
&=\nabla H^*(x(t))\cdot(v^*-x(t))-(q(t)-q^*)\cdot v^*+(\nabla H(x(t))-q(t))\cdot(v(t)-v^*)\hspace{.5cm}\label{q30}\\
&=B_1(t)+B_2(t)+B_3(t),
\end{eqnarray}
where $B_i(t)$, $i\in\left\{1,2,3\right\}$ is the $i$th term in the RHS of (\ref{q30}). Since $x(t)\in V$ and $v^*$ is maximizing $H^*(\cdot)$ over the compact set $V$, then we have 
\begin{equation}
B_{1}(t)\geq H^*(v^*) - H^*(x(t))\geq 0.
\end{equation}
Thus, for any $\epsilon_1>0$, there exist sufficiently small $\epsilon_2>0$ such that
\begin{equation}
\label{eq-b1}
B_{1}(t)\geq \epsilon_2\text{  as long as  }\left\|x(t)-v^*\right\|\geq \epsilon_1.
\end{equation}
Moreover, using complementary slackness \eqn{eq-comp-slack},
\begin{equation}
B_2(t)=-(q(t)-q^*)\cdot v^*=-q(t)\cdot v^*=- \sum_{n\in\mathcal{N}^c}q_n(t)\, v_n^*  \geq 0, \label{aaa38}
\end{equation}
and
\begin{equation}
B_3(t)=\left(\nabla H(x(t))-q(t)\right)\cdot(v(t)-v^*)\geq 0,
\end{equation}
because $v(t)$ maximizes $\nabla H(x(t))-q(t))\cdot v$ over all $v\in V$. 

Now, $\left\|x(t)- v^*\right\|$ must converge to zero (which proves (\ref{eq:converge})).
Indeed, suppose not. Then, there exists $\epsilon_1>0$ and a sequence $t_n, ~n=1,2,\ldots$, $t_n \uparrow \infty$,
such that 
$\left\|x(t_n)-v^*\right\|\geq 2\epsilon_1$. Since $x(t)$ is Lipschitz continuous, this implies that for some $\delta>0$,
$$
\left\|x(t)-v^*\right\|\geq \epsilon_1, ~~ t_n \le t \le t_n+\delta, ~~\forall n,
$$
and then, by \eqn{eq-b1}, for some $\epsilon_2>0$,
$$
B_{1}(t)\geq \epsilon_2, ~~ t_n \le t \le t_n+\delta, ~~\forall n.
$$
This means  that $\int_0^\infty(d\slash dt)F^*(x(t),q(t))=\infty$ (recall that are non-negative), and therefore
$F^*(x(t),q(t)) \to \infty$.
But, this is impossible, because, by the definition of function $F^*$ and Lemma~\ref{lemma3}, 
$\sup_{t\ge 0} \|F^*(x(t),q(t))\| < \infty$. The contradiction proves (\ref{eq:converge}).
\iffalse
Non-negativity of $B_1(\cdot)$, $B_2(\cdot)$ and $B_3(\cdot)$ along with Lipschitz continuity of $x(t)$ show that $\left\|x(t)- v^*\right\|$ must converge to zero, because otherwise $\int_0^\infty(d\slash dt)F^*(x(t),q(t))=\infty$. (This is impossible, since $F^*(x(t),q(t))$ is a uniformly bounded function.) This proves (\ref{eq:converge}).
\fi
\endproof

\subsection{Mapping of the Virtual Matching System of Section~\ref{secmatchvirtualdef} into EGPD Framework}\label{secmatchmodelenhance}

Now we are in position to show that Algorithm~\ref{alg-virtual} for the control of the virtual system in the original matching model in Section~\ref{secmatchmodel} is a special case of EGPD Algorithm.
The mapping of the virtual system of Section~\ref{secmatchvirtualdef} into 
the more general model of Section~\ref{secThemodel} is as follows. 
Consider the following system, which we refer to as a modification of the virtual system.
Suppose the item types $\mathcal{I}$ are modelled as free nodes and let the set of matchings $\mathcal{J}$ be the set of controls $K$. Let us add one \emph{constrained} node per each matching $j\in \mathcal{J}$. (These additional nodes are the \emph{utility} nodes in the terminology of GPD algorithm \cite{stolyar2005maximizing}.) From this point on, for convenience of the notations, we replace the set of indices of item types $\mathcal{I}$ with $\{J+1, \cdots J+I\}$ and denote by $\mathcal{I}^c = \{1,\cdots,J\}$ the set of all constrained nodes. For the constrained nodes we adopt the convention that they never receive any inputs,
i.e. $\lambda_j(t) \equiv 0, ~ j \in \mathcal{I}^c$. We also fix a sufficiently large $c>0$, so that 
$w_j - c < 0$ for all constrained nodes, 
and for each constrained node (or, matching) $j \in \mathcal{I}^c$ we set by convention 
$\mu_j(j) = c - w_j > 0$ and $\mu_i(j) = c > 0, ~i \in \mathcal{I}^c\setminus \{j\}$.
These conventions about the constrained nodes guarantee that under any control strategy, their
queues are automatically stable. In fact, for any $i\in\mathcal{I}^c$ and any initial value $ Q_i(0)$, the queue length  $Q_i(t)$ will decrease until it hits 0 within a finite time and then it will remain at 0. This allows 
to assume, without loss of generality,  that $Q_i(t) \equiv 0$ for all constrained nodes . 

For a matching (or, constrained node) $j$, we have $b(j,t) = (\lambda_i(t)-\mu_i(j),~i\in\mathcal{I}^c \cup \mathcal{I})$ and $\bar{b}(j)=\mathbb{E} b(j,t)$. Note that, for $i\in\mathcal{I}^c$, 
$\bar{b}_i(j)= w_j - c$ if $i=j$ and $\bar{b}_i(j) =-c$ otherwise. 

If $j(t)$ is the matching chosen at $t$, then the compact rate region $V\subset \mathbb{R}^{J+I}$ is the set of all possible vectors
$$
v=(v_1, \ldots, v_J, v_{J+1}, \ldots, v_{J+I})
$$
being possible long-term average values of $\bar{b}(j(t))$ under different matching strategies
(see formal definition in Section \ref{rateregion}). 

Finally, we define the utility function $H(v)$ as follows:
$$
H(v_1, \ldots, v_J, v_{J+1}, \ldots, v_{J+I}) = G(v_1+c, \ldots, v_J+c).
$$

Given these conventions, it is easy to see that the problem of maximizing $G(X_1, \ldots, X_J)$
(subject to the stability of the queues) in the original matching system is equivalent to the problem
of maximizing $H(X_1, \ldots, X_J, X_{J+1}, \ldots, X_{J+I})$ 
(subject to the stability of the queues) 
in the modified system defined in this subsection. The latter system is a special case of the general system
of Section~\ref{secThemodel}. If we specialize Algorithm~\ref{alg-egpd} to the modified system, and then rewrite it in terms of
the original virtual system, we obtain Algorithm~\ref{alg-virtual}.
 Assumption~\ref{assumption 1}, specialized to the modified system and expressed in terms 
of the original virtual system, gives the formal meaning of 
Assumption~\ref{assumption-matching} (which is stated informally).

The mapping described in this section and the asymptotic optimality of the EGPD algorithm 
under Assumption~\ref{assumption 1} imply the asymptotic optimality of Algorithm~\ref{alg-virtual}
under Assumption~\ref{assumption-matching}.

%%%%%%%%%%%%%%%%%%%%%%%%%%%%%%%%%%%%%%%%%%%%%%%%%%%%%%%%%%%%%%%%%%%%%%%%
\section{Simulations}\label{sec:experiment}

In this section, we evaluate the performance of EGPD algorithm via simulations. 
Consider the system described in Section \ref{sec.introduction}. We extend the set of possible matchings  
by including ``single'' matchings
(see Section~\ref{sec-assump-nonrestrict}):
\begin{equation}
\left\{\langle\varnothing\rangle,\langle1\rangle,\langle2\rangle,\langle3\rangle,\langle4\rangle,\langle1,2\rangle,\langle2,3\rangle,\langle2,3,4\rangle\right\}\nonumber.
\end{equation} 

The reward vector is $w=(0,-1,-1,1,2,5,4,7)$ where its $j$-th component corresponds to the $j$th element of the matchings' set.
We consider a linear utility function, namely the sum of average rewards due to different matchings.   
The vector of arrivals rates is $\alpha=(1.2,1.5, 2,0.8)$.

For our linear utility function, the EGPD algorithm for the virtual system is given by rule \eqn{rule3-linear}.

\begin{figure}[htb]
    \begin{subfigure}[t]{.4\textwidth}
        \centering
        \caption{Virtual System}
		\includegraphics[width=.97\textwidth]{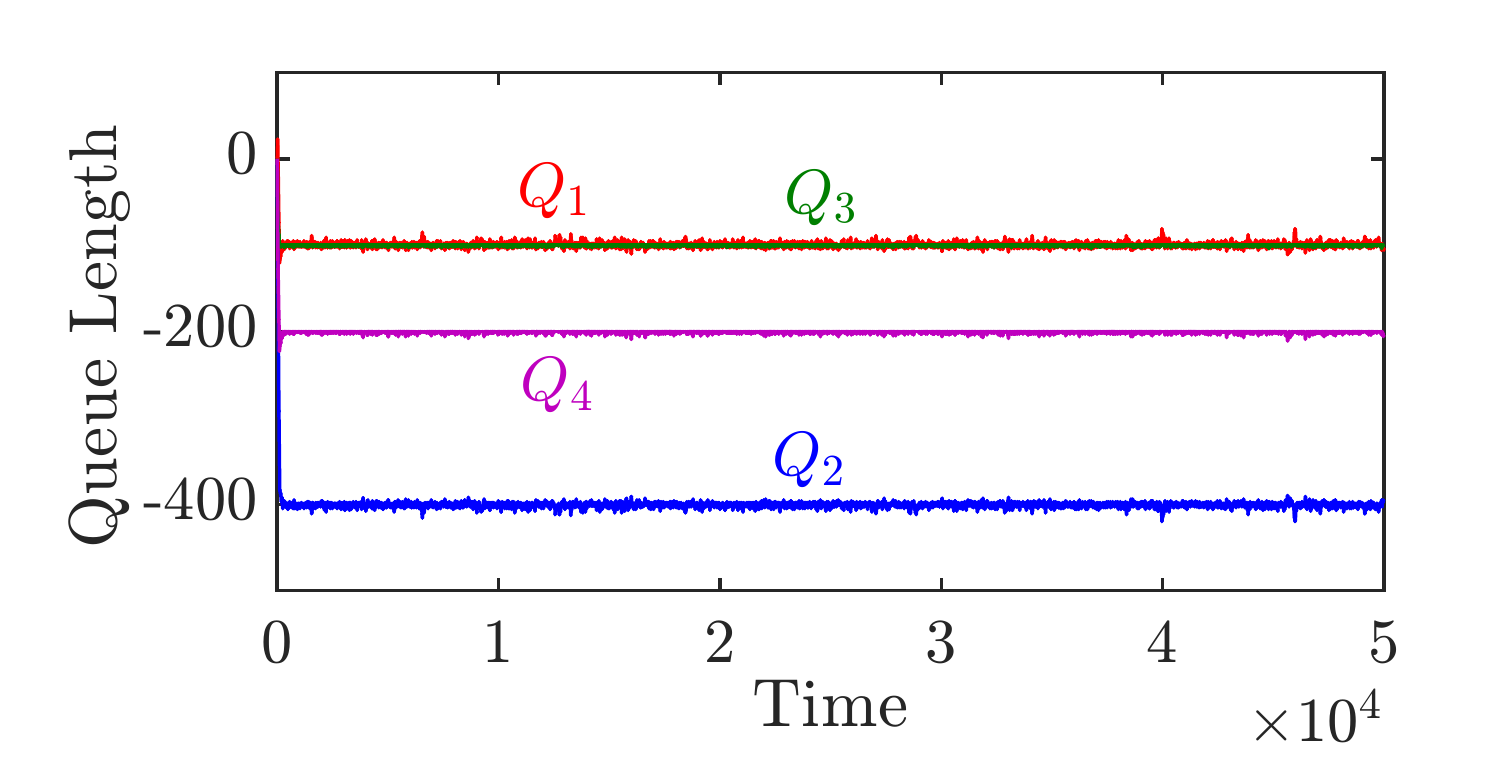}
    \end{subfigure}%
        \centering
    \centering
    \begin{subfigure}[t]{.4\textwidth}
        \centering
        \caption{Physical System}
		\includegraphics[width=.97\textwidth]{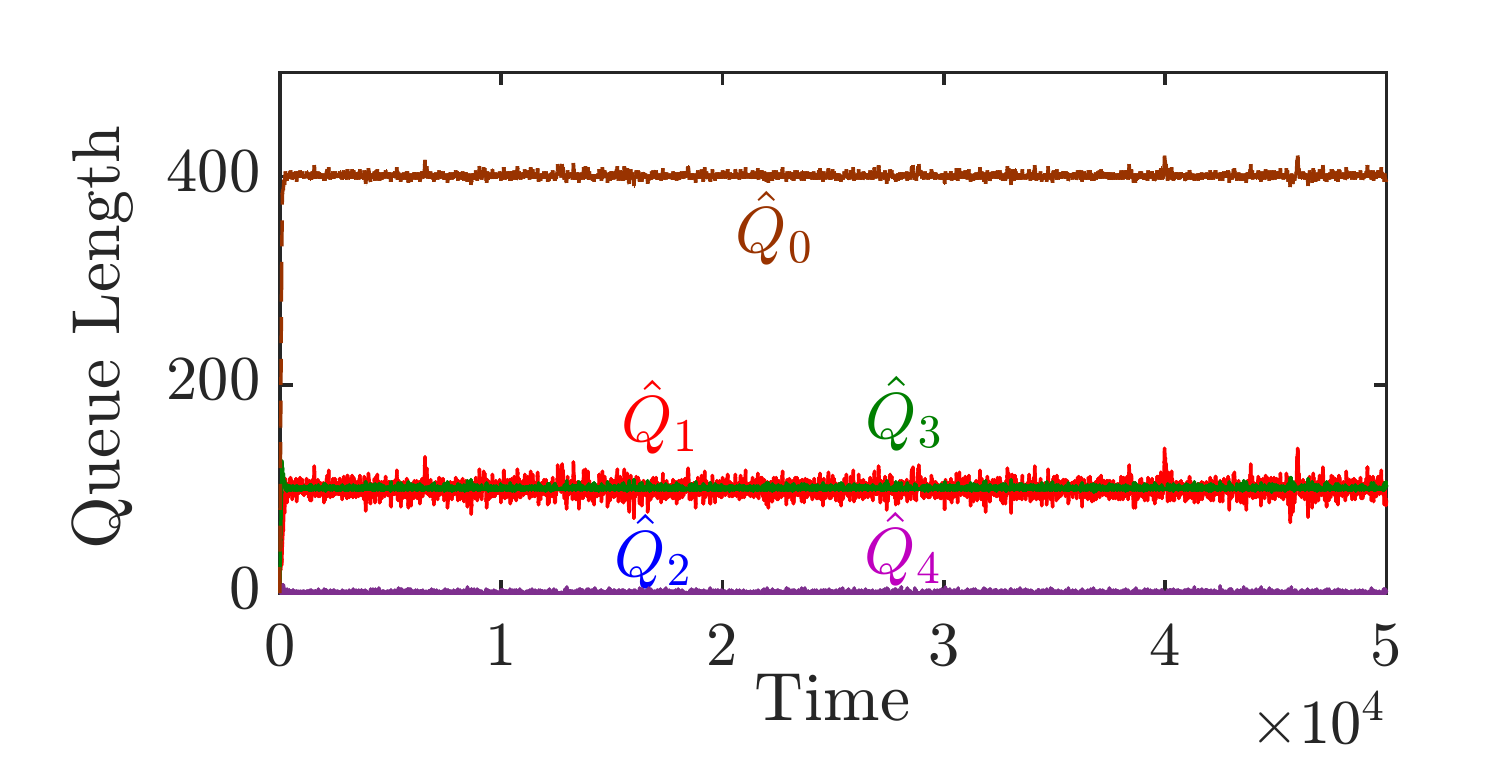}		
    \end{subfigure}  
\caption{Queue trajectories of the virtual and physical systems under EGPD algorithm.}
\label{fig:EGPD_trajectories}
\end{figure}

\textbf{A. Average reward maximization.} We use parameter $\beta=0.01$.
Figure~\ref{fig:EGPD_trajectories} shows the queue trajectories of the virtual and physical systems under the EGPD algorithm. All queues are initially empty. We observe that all queues are quickly ``converging''. Nearly all type 2 and 4 items are matched right after they enter the system, while there exist a queue of around 100 items of types 1 and 3.

The rates at which matchings are activated under EGPD algorithm are provided in table \ref{table4}, which shows that these rates are close to the optimal ones, obtained by solving the underlying optimization problem (which 
is a linear program in this case). 
Therefore, as expected, the algorithm yields near optimal performance for small $\beta$.
Note that solving the optimization problem requires the knowledge 
of arrival rates (as well as other system parameters), while our algorithm need not know arrival rates.

%%%%%%%%%%%%%%%%%%%%%%%%%%%%%%%%%%%%%%%%%%%%%%%%%%%%%%%%%%%%%%%%%%%%%%%%%%
\begin{table}[htb]
  \centering
  \footnotesize
 \caption{Matching rates: Optimal vs. EGPD. (Runtime=30000)}
	  \label{table4}
    \begin{tabular}{lcccccccc}
        \toprule
    \multirow{2}{*}{Method}&\phantom{a}&\multicolumn{7}{c}{Matchings}\\
    \cmidrule{3-9}
       && $\langle1\rangle$     & $\langle2\rangle$       & $\langle3\rangle$      & $\langle4\rangle$       & $\langle1,2\rangle$   & $\langle2,3\rangle$    & $\langle2,3,4\rangle$  \\
    \midrule
      EGPD  & &0 &0&1.69345 &0.4829 &1.1924 &0&0.31075\\
  Optimal &   &0&0&1.70005 &0.49995 &1.2001&0&0.29975\\
    \bottomrule
    \end{tabular}%
\end{table}%
%%%%%%%%%%%%%%%%%%%%%%%%%%%%%%%%%%%%%%%%%%%%%%%%%%%%%%%%%%%%%%%%%%%%%%%%%%
\begin{figure}[htb]
\centering
	\includegraphics[width=0.35\textwidth]{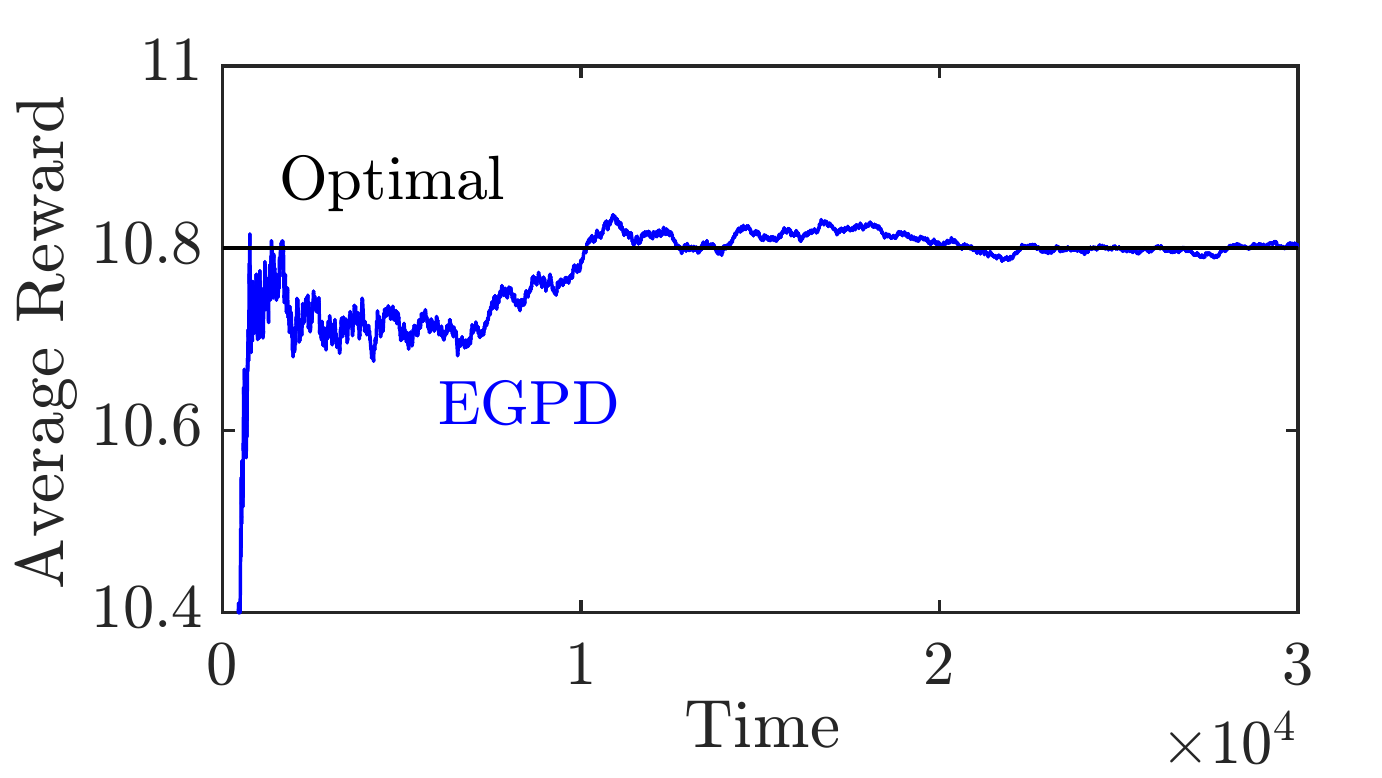}
\caption{Average matching reward under the EGPD algorithm.}
\label{fig:opt}
\end{figure}

Figure~\ref{fig:opt} demonstrates the average matching reward per unit time. We have calculated the optimal average reward (by solving the linear program) which is equal to 10.8, and plotted it on the figure. As clear from the graph, the running average reward under EGPD algorithm is getting very close to optimal objective value and this convergence is sufficiently fast.

\textbf{B. Effect of parameter $\beta$.} In order for $\beta {Q}$ in the virtual system to ``stay close'' to some $q^*\in Q^*$, parameter $\beta$ should be small. Therefore, {\em as long as parameter $\beta$ is sufficiently small}, 
the algorithm is nearly optimal and the virtual queue lengths are roughly of the order $1/\beta$. 
As $\beta$ is increasing, the accuracy of the algorithm in terms of average reward maximization decreases, while the queues become smaller. 

\begin{figure}[htb]
\centering
		\includegraphics[width=0.35\textwidth]{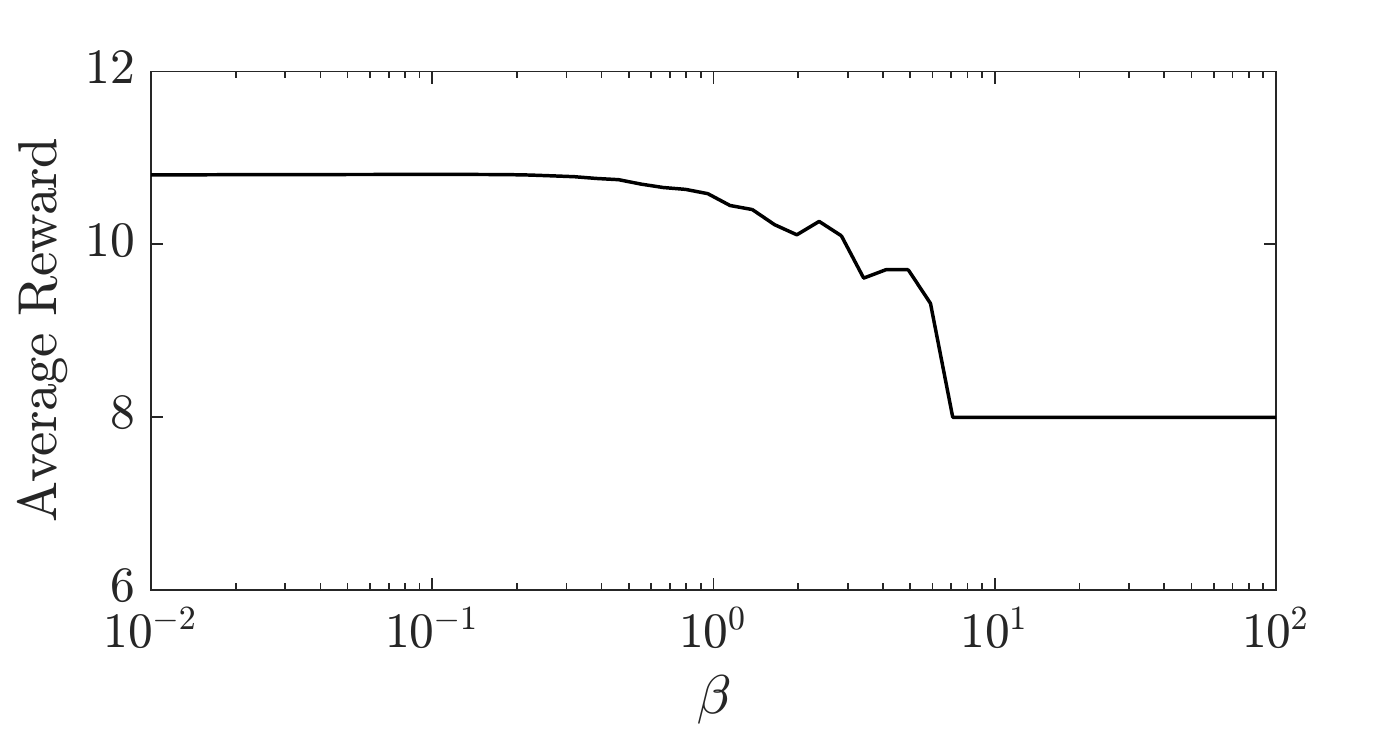}
\caption{Average matching reward for different values of $\beta$.}\label{fig:optbeta}
\vspace{-.3cm}
\end{figure}

The dependence of the average reward on $\beta$ for the considered scenario is 
shown on Figure~\ref{fig:optbeta}. First, we note that the average reward remains nearly optimal for values of $\beta$ almost as large as $1$ (i.e. not even very small in absolute terms). Then, as $\beta$ changes from $1$ to about $10$, the average reward decreases and reaches the lower ``plateau,'' and then remains constant for $\beta\ge 10$. Thus, as expected, the algorithm is effective in terms of reward maximization when $\beta$ is sufficiently small (less than $1$ in our scenario); when $\beta$ is sufficiently large (greater than $10$ in our scenario), the average reward is also roughly independent of $\beta$, but is at a lower, suboptimal level.

We note that larger values of $\beta$ have the benefit of reducing the queues and, as a result, reducing (as we will see next) the algorithm response (or, adaptation) time to changes of the items' arrival rates. (Shorter queue also mean lower holding costs, if such are a part of the model. This will be discussed in Section~\ref{profit}.)
Therefore, the value of parameter $\beta$ should be chosen, very informally speaking, ``as large as possible, but not larger''. 

\begin{figure}[htb]
\centering
    {\begin{subfigure}[t]{.4\textwidth}
        \centering
        \caption{Virtual System}
        \includegraphics[width=.95\textwidth]{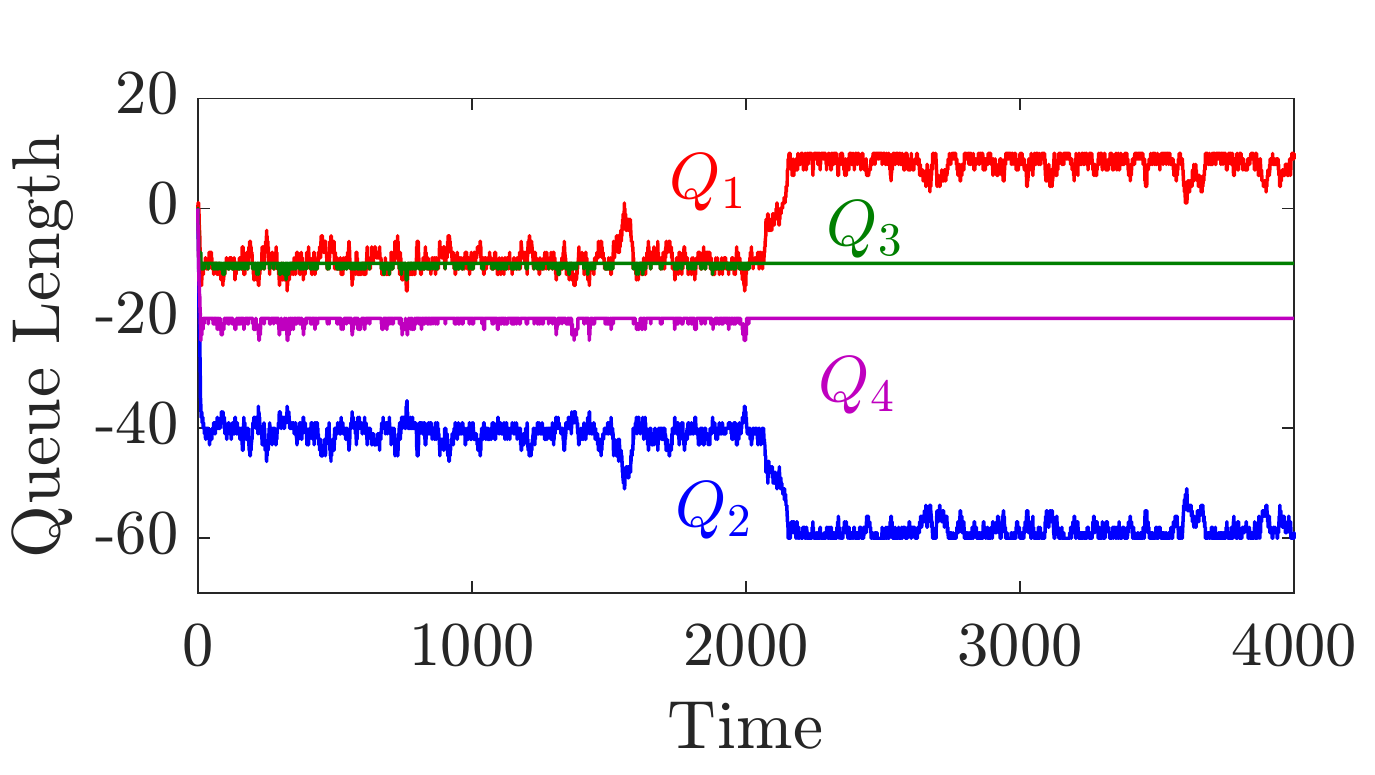}
    \end{subfigure}%
        \begin{subfigure}[t]{.4\textwidth}
        \centering
        \caption{Physical System}
		\includegraphics[width=.95\textwidth]{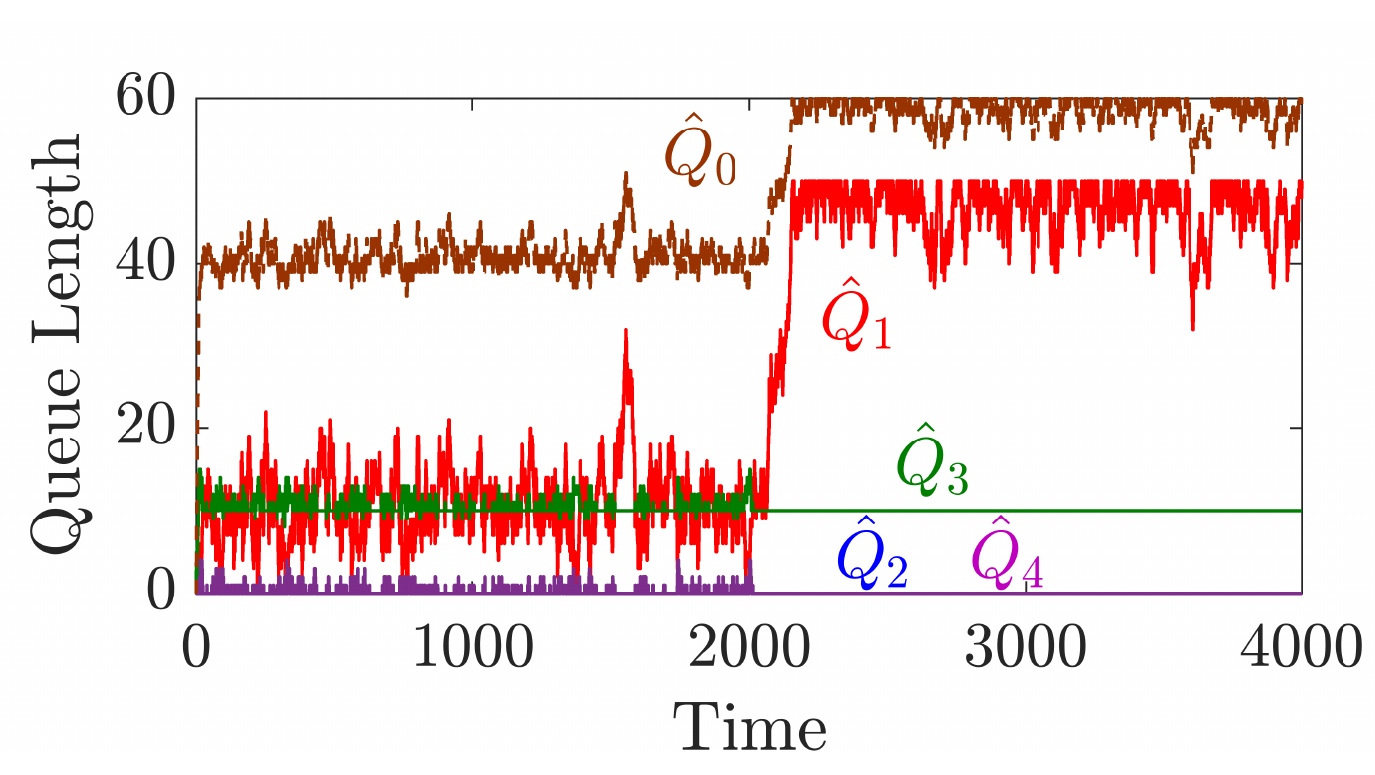}
    \end{subfigure}}%
\caption{Adaptation to the changes in arrival rates.} \label{fig:arrival_change}
\end{figure}

\textbf{C. Automatic adaptation to changes in arrival process.} An important robustness issue is how quickly the EGPD algorithm responds to the changes in the arrival process. In the following experiment, the arrival rates are changed to  $\alpha = (1.8,0.8,1.4,1)$ at time 2000. This change leads to different optimal matching rates and thus different optimal value. If quick response to arrival rate changes is important, a larger $\beta$ is preferable.
Here we use $\beta=0.1$. Figure~\ref{fig:arrival_change} shows the queue trajectories of the virtual and physical systems. We observe that EGPD automatically adapts to the new arrival rates and reaches the new ``right'' queue lengths, without using any a priori information on this change.

\iffalse
\textbf{C. Effect of parameter $\beta$ on average profit.}
Figure~\ref{fig:beta} shows the dependence of average profit on parameter $\beta$. As $\beta$ increases, the average reward is getting far from the optimal one and the average holding cost decreases. We observe that the average profit increases and then decreases for larger values of $\beta$. We conjecture a similar behavior to hold in general. With this regard, we observe that by fixing appropriate $\beta$ (for example, it is around 2 in this example), the average profit obtained from the EGPD algorithm will be increased.
\begin{figure}[h]
		{\includegraphics[width=0.5\textwidth]{change_beta.pdf}}
{Average profit for different values of $\beta$.\label{fig:beta}}{}

\end{figure}
\fi

\section{Heuristics for the Objective Including both Matching Rewards and Holding Costs}\label{profit}

\subsection{General discussion}

The scheme we proposed in Section~\ref{secdynvirtual} for the matching model is asymptotically optimal for the reward maximization problem. (We will refer to this entire scheme as EGPD, because EGPD is its key part,
applied to the virtual system and determining matching choices.)
In practical systems, the objective may be more general, namely maximizing the average ``profit'' defined as average reward minus average queue holding cost. We now informally discuss how EGPD can be used to achieve better profit in the system (even though it is not specifically designed for that).

For the purposes of the discussion below, we assume linear holding costs with rate vector $c = (c_i,\,i\in\mathcal{I})$; that is the average holding cost over interval $[0,T]$ is
\begin{eqnarray}
\frac{1}{T}\int_{0}^T c \cdot \hat Q(t)dt.
\end{eqnarray}
Suppose the arrival rates are scaled up by a factor $r> 0$. This simply speeds up the process $r$ times, so that the average reward increases $r$ times, while the holding cost remains same. Thus for systems with ``high'' arrival rates, the rewards dominate the profit objective and we expect the average profit obtained by the EGPD algorithm to be ``close'' to the optimal one. In other cases, holding costs may dominate, for example when the system is in (appropriately defined) heavy traffic (see \cite{gurvich2014dynamic,busic2014optimization}) -- this makes the queues necessarily large.
When the optimal average rewards and optimal holding cost are on the same scale, the EGPD parameter settings can be used to control the tradeoff between these two performance measures, thus potentially improving the average profit. 
We now briefly discuss different heuristic approaches for profit improvement within the framework of our scheme.

\textbf{Choice of parameter $\beta$.} As discussed in Section~\ref{sec:experiment}, 
as long as parameter $\beta$ is sufficiently small, 
the virtual queue lengths under EGPD are large, roughly of the order $1/\beta$. To see how this affects the holding cost, consider two cases:
 \begin{enumerate}[label=(\roman*),font=\itshape]
 \item If ${Q}_i(t) \geq 0$, then $\hat Q_i(t)$ will also be large (of the order of at least $1/\beta$) since the inequality $\hat Q_i(t)\geq {Q}_i(t)$ holds for all $i\in\mathcal{I}$ at all $t$.
 \item If ${Q}_i(t) < 0$, this has an indirect impact on the holding cost. In particular, large $|{Q}_i(t)|$ in this case would imply more incomplete matchings. This subsequently results in a higher holding cost.
\end{enumerate} 	  
Therefore, parameter $\beta$ should be chosen as large as possible, but not to exceed the level beyond which the average rewards start to be significantly (negatively) affected.

\textbf{Additional queue scaling.} Consider arbitrary positive weights $\gamma_i, ~i\in \mathcal{I}$. All the results for the EGPD algorithm hold if we use more general rule 
\begin{equation}
j(t)\in\argmax_{j\in \mathcal{J}}  \left[ (\partial G(X(t))/\partial X_j)\, w_j+\sum_{i\in{\mathcal{I}}}\beta\, \gamma_i \, {Q}_i(t)\, \mu_i(j)\right].\label{weighted}
\end{equation}
instead of (\ref{rule3}). In this case, it is the weighted vector $(\gamma_i \beta {Q}_i(t),\, i\in \mathcal{I})$ 
(not $\beta {Q}(t)$) that will be close to an optimal dual solution $q^*$. This property may be used to reduce the holding cost by giving higher weights to more ``expensive'' queues (with large $c_i$), thus making them relatively smaller.

\textbf{Matching completion order.} There is a flexibility in choosing which incomplete matching to complete first. For the average matching 
reward maximization this does not matter (so, earlier we specified FCFS rule for concreteness).
However, if the holding costs are a consideration, 
one may pick incomplete matchings with higher associated holding cost to be completed first. 

\subsection{Simulation: Average Profit in a Bipartite Matching System}
\label{subsec-bipart-sim}

\begin{figure}[h]
\centering
		\includegraphics[width=0.6\textwidth]{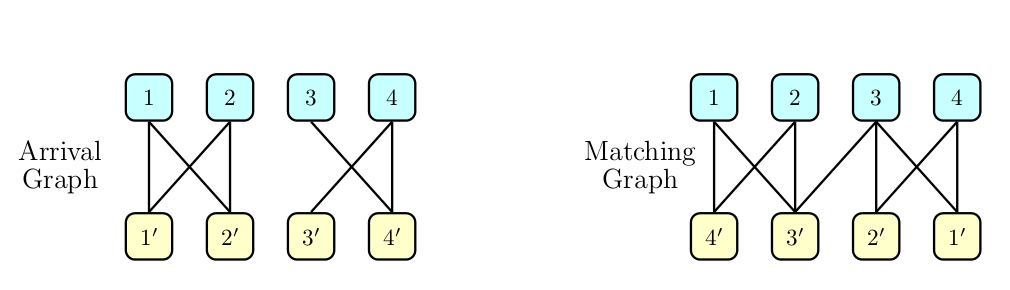}	
\caption{Illustration of the matching system.}
\label{figmatch2}
\end{figure}

Consider a bipartite matching system, where items arrive in pairs, and the matchings are pairs as well.
It is depicted in 
Figure~\ref{figmatch2}. There are 8 item types $\{1,2,3,4,1',2',3',4'\}.$ The arrival graph is on the left, where each edge shows a possible arrival pair, and the plot in the right hand side is the matching graph with edges representing the possible matchings. Up to two matchings can be done per each arrival ($m=2$).

We consider the process in discrete time $t=1,2,\cdots$. The arrival process is i.i.d. across time.
Specifically, at each time $t$, a pair of items enters the system. The probabilities (rates)
of different arrival pairs are specified in Table~\ref{table5}.

\begin{table}[H]
\footnotesize
  \centering
    \caption{Probabilities (rates) of different arrival pairs.}
  \label{table5}
  {\renewcommand{\arraystretch}{1.1}%
    \begin{tabular}{c|ccccccc}
    \hline
{Arrival pairs} & (1,$1'$) & (1,$2'$) & (2,$1'$) & (2,$2'$) & (3,$4'$) & (4,$3'$) & (4,$4'$) \\ \hline
{Probability} & 0.166 & 0.083 & 0.087 & 0.083 & 0.2324 & 0.2656 & 0.083\\
 \hline
\end{tabular}}
\end{table}%
\vspace{-.2cm}

It is easy to check that this system satisfies necessary and sufficient condition 
\cite{busic2010stability} for bi-partite matching systems
to be stabilizable. The condition, called NCond in \cite{busic2010stability}, is as follows. Suppose the matching graph is connected. 
Consider a subset $T$ of items from the top part of the bi-partite graph, and denote by $\alpha_T$ the total arrival rate of all items in $T$. Denote by $B(T)$ the subset of items from the bottom part of the graph that can be matched with at least one item in $T$, and by $\alpha_{B(T)}$ the total arrival rate of all items in $B(T)$.
Then, the system is stabilizable if and only if $\alpha_T < \alpha_{B(T)}$ holds for any {\em strict} subset $T$ of ``top'' items. 

Now, let us see if our Assumption~\ref{assumption-matching} holds.
Since this is a bipartite matching system, with items arriving and departing in pairs, virtual queues satisfy the following linear relation $Q_1(t)+Q_2(t)+Q_3(t)+Q_4(t) - (Q_{1'}(t)+Q_{2'}(t)+Q_{3'}(t)+Q_{4'}(t)) \equiv 0$.
However, given that NCond condition holds,
it is easy to see that our Assumption~\ref{assumption-matching} (formally given by 
Assumption \ref{assumption 1}) holds for this system in the sense 
described in Section~\ref{sec-assump-nonrestrict}, namely after an orthogonal change of coordinates.
(We emphasize again that the algorithm itself remains as is, it does {\em not} need to do any change 
of coordinates.) Therefore, the EGPD algorithm is asymptotically optimal for this system 
for the average reward maximization objective.

Assume linear holding costs, $c\cdot \hat Q(t)$, with the cost rate vector $c=(0.1,0.2,0.3,0.4,0.4,0.3,0.2,0.1)$. 
The matching rewards for different matchings are given in Table~\ref{table:matchpair}.

\begin{table}[H]
\footnotesize
  \centering
    \caption{Matching rewards. }
  \label{table:matchpair}
  {\renewcommand{\arraystretch}{1.1}%
  \setlength\tabcolsep{4pt}
    \begin{tabular}{c|ccccccccc}
\hline        
{Matchings} & $\langle 1,3'\rangle$ & $\langle 1,4'\rangle$ & $\langle 2,3'\rangle$ & $\langle 2,4'\rangle$ & $\langle 3,1'\rangle$ & $\langle 3,2'\rangle$ & $\langle 3,3'\rangle$ & $\langle 4,1'\rangle$ & $\langle 4,2'\rangle$ \\ \hline
{Reward} & 5 & 50 & 5 & 50 & 5 & 50 & 5& 50&5 \\
\hline
\end{tabular}}
\end{table}%

\begin{figure}[!h]
\centering

    \begin{subfigure}[b]{0.6\textwidth}
        \centering
        \includegraphics[width=.65\textwidth]{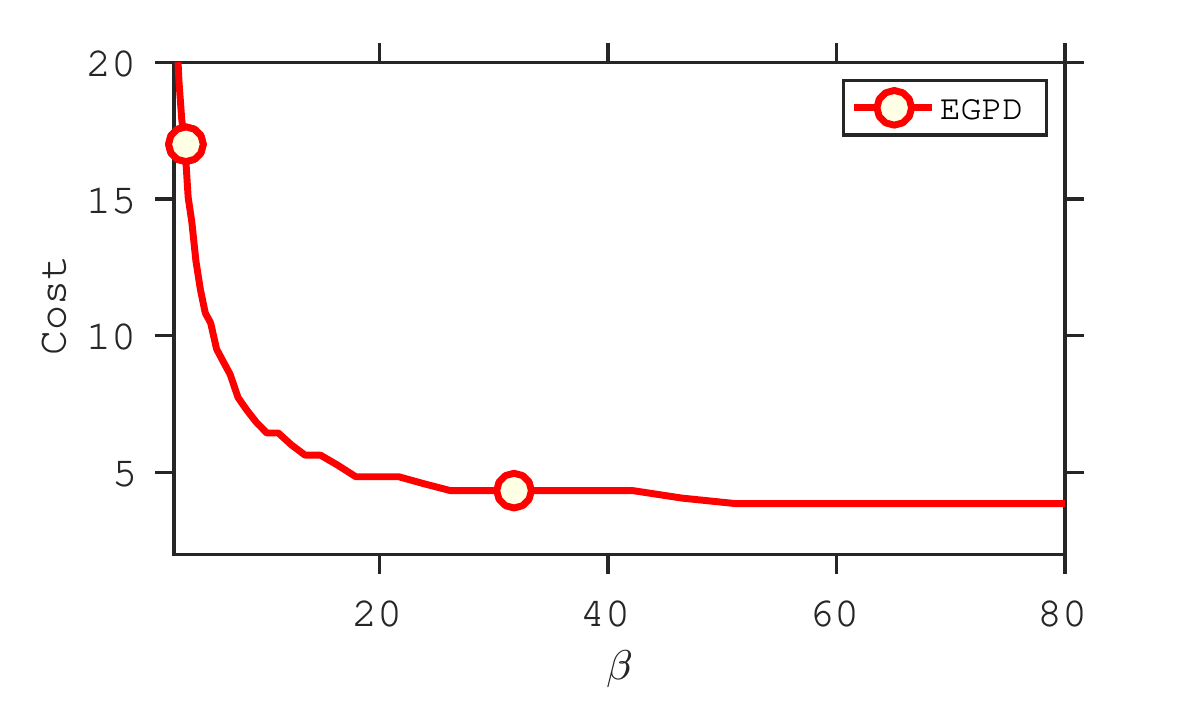}
        \caption{Average holding cost}
    \end{subfigure} 
        
    \begin{subfigure}[b]{0.6\textwidth}
        \centering
        \includegraphics[width=.65\textwidth]{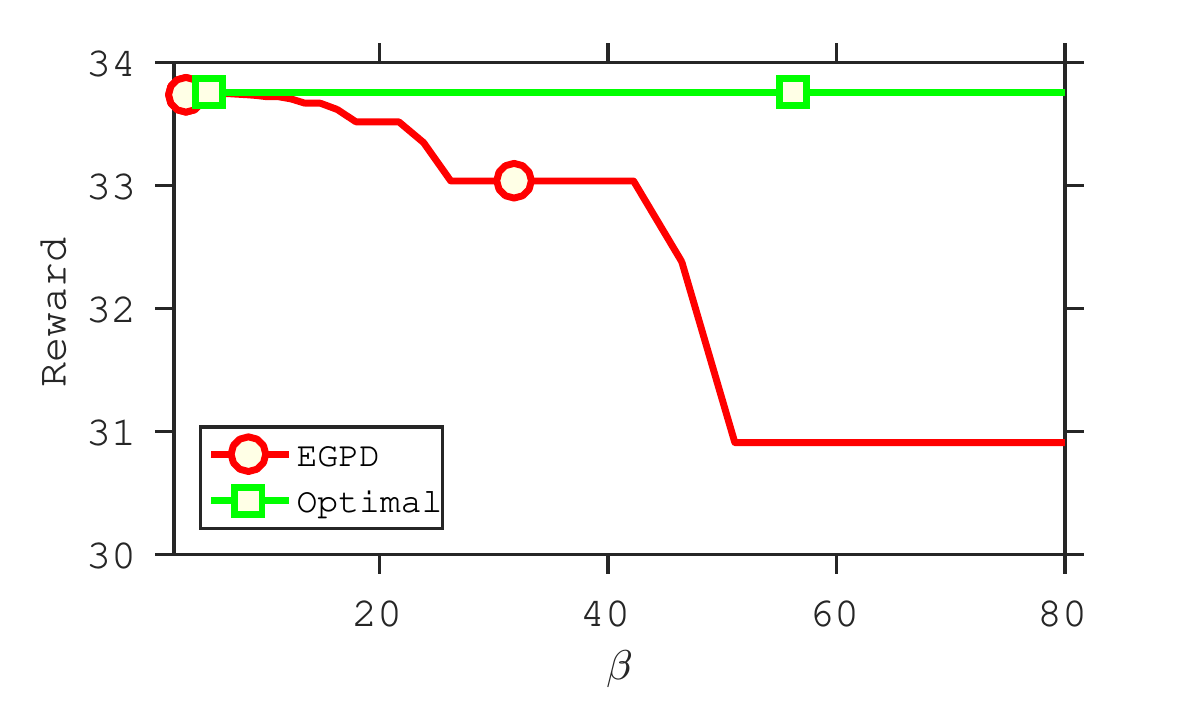}
        \caption{Average reward}
    \end{subfigure}   
    
    \begin{subfigure}[t]{0.6\textwidth}
        \centering
        \includegraphics[width=.65\textwidth]{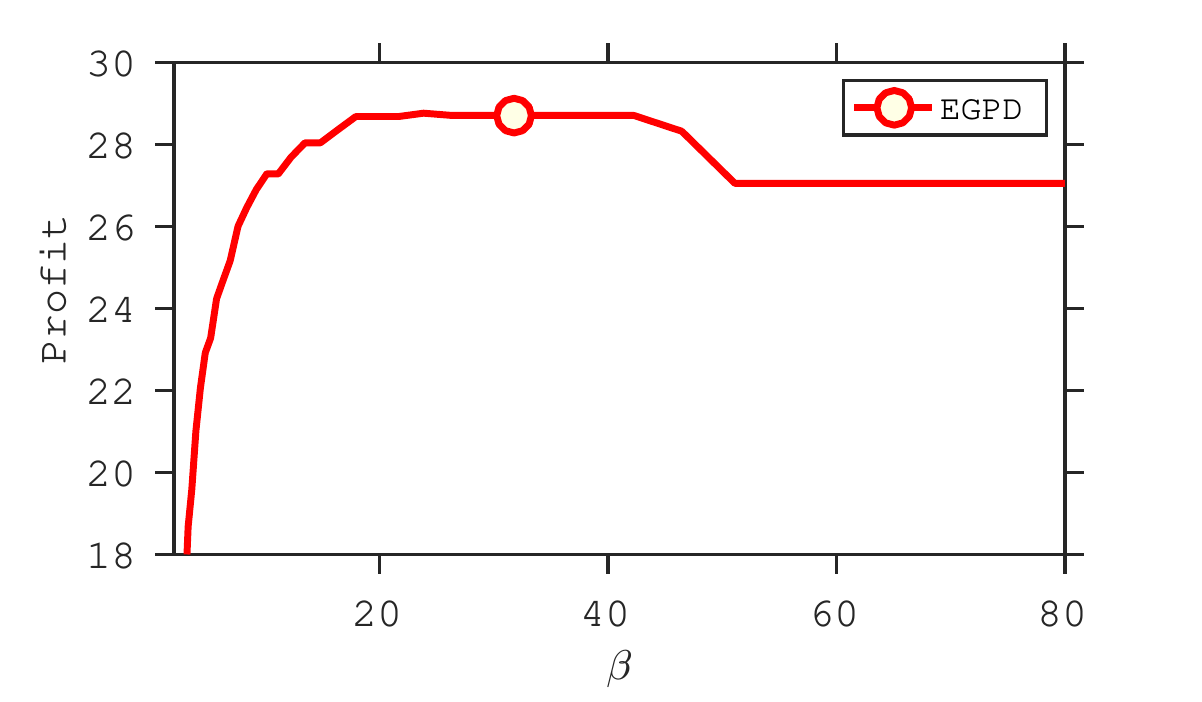}        		     

        \caption{Average profit}
    \end{subfigure}
    \caption{EGPD algorithm performance.}
    \label{comparison}
\end{figure}

We simulated this system under EGDP scheme.
Figure~\ref{comparison} shows the dependence of EGPD average performance metrics on the parameter $\beta$. The range of $\beta$ is shown within which the average reward declines from its optimal (largest) value to the ``plateau'' it reaches when $\beta$ is large. Parts (a), (b) and (c) show average holding cost, reward and profit, respectively; the average profit is the average reward minus the average holding cost.

We see that the average profit is maximized within a certain range of values of $\beta$, where, roughly speaking, the average reward is ``still'' close to optimal and the average holding cost is ``already'' close to the best achievable by EGPD. 
We conjecture that the average profit with such choice of $\beta$ is reasonably close to the optimal profit under any control algorithm. Verifying and quantifying this informal conjecture is an interesting subject for future research.

%**************
% Conclusion
%**************
\section{Conclusions}
\label{sec-conc}

In this paper we have proposed an approach for optimal dynamic control of general matching systems. The central idea is using a virtual matching system allowing negative (as well as positive) queues, as part of the overall control scheme.  The virtual system fits into a queueing network framework, except the queues may be negative, and it is controlled by 
an extended version of the GPD algorithm, called EGPD. We prove EGPD asymptotic optimality.
The approach is very generic, not restricted to special cases, such as bipartite matching. 
The proposed scheme is also very robust in the sense that it does not require the knowledge of input rates, and automatically adapts to changing input rates. Simulations demonstrate good performance of the algorithm.

Although the scheme that we develop has the average reward maximization as its objective, the parameter setting can be used to achieve good performance in terms of the more general objective, which includes holding costs. 
Addressing this and other more general objectives within a dynamic control framework, not requiring a priori knowledge of the item arrival rates, is an important future subject. 

%**************
% Bibliography
%**************
%\nocite{*}
%\newpage
\bibliographystyle{plainnat}
%\bibliographystyle{siam}
%\bibliography{research}

\begin{thebibliography}{}

\bibitem{adan2015reversibility}
Adan, I., Bu{\v s}i\'c, A., Mairesse, J., Weiss, G.: Reversibility and further
  properties of FCFS  infinite bipartite matching.
\newblock Mathematics of Operations Research  (2015)

\bibitem{adan2012exact}
Adan, I., Weiss, G.: Exact FCFS matching rates for two infinite multitype
  sequences.
\newblock Operations research \textbf{60}(2), 475--489 (2012)

\bibitem{buke2015stabilizing}
B{\"u}ke, B., Chen, H.: Stabilizing policies for probabilistic matching
  systems.
\newblock Queueing Systems \textbf{80}(1-2), 35--69 (2015)

\bibitem{busic2010stability}
Bu{\v s}i\'c, A., Gupta, V., Mairesse, J.: Stability of the bipartite matching
  model.
\newblock Advances in Applied Probability \textbf{45}(2), 351--378
  (2013)

\bibitem{busic2014optimization}
Bu{\v s}i\'c, A., Meyn, S.: Optimization of dynamic matching models.
\newblock arXiv preprint arXiv:1411.1044  (2014)

\bibitem{rene2009fcfs}
Caldentey, R., Kaplan, E.H., Weiss, G.: FCFS infinite bipartite matching of
  servers and customers.
\newblock Advances in Applied Probability \textbf{41}(3), 695--730 (2009)

\bibitem{gurvich2014dynamic}
Gurvich, I., Ward, A.: On the dynamic control of matching queues.
\newblock Stochastic Systems \textbf{4}(2), 479--523 (2014)

\bibitem{kashyap1966double}
Kashyap, B.: The double-ended queue with bulk service and limited waiting
  space.
\newblock Operations Research \textbf{14}(5), 822--834 (1966)

\bibitem{lovasz2009matching}
Lov{\'a}sz, L., Plummer, M.D.: Matching theory, vol. 367.
\newblock American Mathematical Soc. (2009)

\bibitem{mairesse2014stability}
Mairesse, J., Moyal, P.: Stability of the stochastic matching model.
\newblock Journal of Applied Probability \textbf{53}(4), 1064--1077 (2016)

\bibitem{mehta2012online}
Mehta, A.: Online matching and ad allocation.
\newblock Theoretical Computer Science \textbf{8}(4), 265--368 (2012)

\bibitem{PW2008}
Plambeck, E.L., Ward, A.R.: Optimal control of a high-volume assemble-to-order
  system with maximum leadtime quotation and expediting.
\newblock Queueing Systems \textbf{60}(1), 1--69 (2008)

\bibitem{stolyar2005maximizing}
Stolyar, A.L.: Maximizing queueing network utility subject to stability: Greedy
  primal-dual algorithm.
\newblock Queueing Systems \textbf{50}(4), 401--457 (2005)

\bibitem{St2005gpdgen}
Stolyar, A.L.: Greedy primal-dual algorithm for dynamic resource allocation in
  complex networks.
\newblock Queueing Systems \textbf{54}(3), 203--220 (2006)

\bibitem{stolyar2010control}
Stolyar, A.L., Tezcan, T.: Control of systems with flexible multi-server pools:
  a shadow routing approach.
\newblock Queueing Systems \textbf{66}(1), 1--51 (2010)

\end{thebibliography}

\end{document}